\documentclass[draftclsnofoot,onecolumn, 12pt]{IEEEtran}

\usepackage{amsmath,bm}
\usepackage{mathtools}
\mathtoolsset{showonlyrefs,showmanualtags}
\usepackage{graphicx}
\usepackage{amssymb}
\usepackage{bm}
\usepackage{cite}
\usepackage[thinc]{esdiff}
\usepackage[ruled]{algorithm2e}
\usepackage{xcolor}
\newtheorem{example}{Example}
\newtheorem{remark}{Remark}
\newtheorem{theorem}{Theorem}
\newtheorem{lemma}{Lemma}
\newtheorem{cor}{Corollary}

\newcommand{\I}{\mathrm{I}}
\newcommand{\E}{\mathrm{e}}

\begin{document}

\title{A Method For Bounding Tail Probabilities}

\author{Nikola Zlatanov 
\thanks{This paper is published in IEEE Access, see https://doi.org/10.1109/ACCESS.2026.3650974}
\thanks{N. Zlatanov is with the Faculty of Computer Science and Software Engineering, Innopolis University, Innopolis, Russia. E-mails: n.zlatanov@innopolis.ru}
 }

\maketitle

\begin{abstract}
We present a  method for upper and lower bounding the right and the left  tail probabilities of continuous random variables (RVs). For the  right tail probability of RV $X$ with probability density function $f (x)$, this method requires first setting a continuous, positive, and strictly decreasing function $g (x)$  such that $-f (x)/g' (x)$ is a decreasing and increasing function,  $\forall x>x_0$, which results in  upper and lower bounds, respectively,   given in the form $-f (x) g (x)/g' (x)$, $\forall x>x_0$, where $x_0$ is some point.  Similarly,  for the upper and lower bounds on the left tail probability of $X$, this method requires first setting a continuous, positive, and  strictly increasing function $g (x)$  such that $f (x)/g' (x)$ is an increasing and decreasing function,  $\forall x<x_0$, which results in  upper and lower bounds, respectively,   given in the form $f (x) g (x)/g' (x)$, $\forall x<x_0$.
We provide some examples of good candidates for the function $g (x)$. 
We also establish connections between the new bounds and Markov's inequality and Chernoff's bound. In addition, we provide an iterative method for obtaining ever tighter lower and upper bounds, under certain conditions.   As an application, we use the proposed method to derive a novel closed-form asymptotic expression of the converse bound  on the capacity of the additive white Gaussian noise (AWGN) channel in the finite-blocklength regime, which is tighter than the closed-form  asymptotic expression by Polyanskiy-Poor-Verdú.
Finally, we provide numerical examples where we show the tightness of the bounds obtained by the proposed method. 
\end{abstract}


\section{Introduction}
 
The most well known and the most utilized methods  for bounding tail probabilities are based on  variations of Markov's inequality \cite{markov_book}.
Markov's inequality   relates the right tail probability of a non-negative random variable (RV) to its mean. The Bienaymé-Chebyshev's inequality  \cite{bienayme1853considerations, chebyshev1867valeurs},   relates the right tail probability of a  RV to its mean and variance, and this inequality can be obtained from Markov's inequality. Other notable  bounds on the tail probabilities that  are based on   Markov's inequality are the Chernoff-Cramér  bound \cite{chernov} and   Hoeffding's inequality \cite{hoeff}, among the most famous. 

Additional tail bounding methods include martingale methods  \cite{mcdiarmid_1989}, information-theoretic methods \cite{ahlswede1976bounds, 10.5555/1146355},  the entropy method  based on logarithmic Sobolev inequalities  \cite{ledoux_1997},  Talagrand’s
induction method \cite{talagrand1995concentration}, etc.
For an  overview of tail bounding  methods, please refer to \cite{con_in_book}. 

Tail bounds are especially important in communications and information theory. For example,   bit/symbol error rates of communications channels corrupted by additive white Gaussian noise (AWGN) are almost always obtained as a function of the Gaussian right tail probability, see \cite[Chapter 4]{proakis2008digital}. On the other hand, the Polyanskiy-Poor-Verdú converse bound for the finite blocklength  AWGN channel capacity, derived in \cite{5452208}, is given in the form of the left tail of the non-central chi-squared distribution, see \cite{7303958, 7589108}. Therefore, having tight bounds on the right and the left tail probabilities would lead to better understanding of results in communications and information theory.

In this paper, we provide a general method for upper and lower bounding both the right tail and the left tail of continuous RVs. In summary, the upper and the lower bounds on the right tail of a continuous RV $X$ with  probability density function (PDF) $f (x)$ and  support on $(l,r)$, are given by
 $$P (x)=-f (x)\frac{g (x)}{g' (x)}, \;\forall x>x_0,$$
 where $g (x)$ is any   continuous, positive, and strictly decreasing function, $\forall x>x_0$, that results in   $\lim\limits_{x\to r}P (x)=0$ being met, when 
$$\frac{f (x)}{-g' (x)}$$
is a decreasing and increasing function $\forall x>x_0$, respectively.  
 
 Similarly, the upper and the lower bounds on the left tail of $X$ are given by
 $$P (x)=f (x)\frac{g (x)}{g' (x)}, \;\forall x<x_0,$$
 where $g (x)$ is any   continuous, positive, and strictly increasing function, $\forall x<x_0$, that results in   $\lim\limits_{x\to l}P (x)=0$ being met, when 
$$\frac{f (x)}{g' (x)}$$
is an increasing and decreasing function $\forall x<x_0$, respectively.
 
 The method is general since there are many functions, $g (x)$, which satisfy the above descriptions and are therefore good candidates for building upper and lower bounds on the right and the left tails. For example, $g (x)=f (x)$ and $g (x)=(x-l)f (x) $ are two possible candidates that lead to tight upper bounds on RVs with exponential and sub-exponential decay of their right tails, respectively. Moreover, $g (x)=(x-l)f (x)$ is a good candidate that leads to a tight upper bound on the left tail. 
 
We also establish  connections between the bounds resulting from the proposed method with Markov's inequality and Chernoff's bound.

 In addition, we provide an iterative method that leads to ever tighter upper and lower bounds on both tails, under certain conditions, using an iterative function of the form
  $$P_{i+1}(x)=\pm f (x)\frac{P_i(x)}{P'_i(x)},$$
for which  the seed is given by  
    $$P_0(x)=\pm f (x)\frac{g (x)}{g' (x)}.$$
    
As an application, we use the proposed tail bounding method to derive tight lower and upper bounds on the converse bound of the AWGN channel capacity in the finite-blocklength regime, derived in \cite{5452208}. It is well known that evaluating the  converse bound of the AWGN channel capacity in the finite-blocklength regime is difficult, and is therefore the main theme  of  papers such as \cite{7303958} and \cite{7589108}. Hence, the tight bounds that we provide are helpful in evaluating this converse bound. Moreover, using the proposed tail bounding method, we derive  a novel closed-form asymptotic expression of the converse bound of the AWGN channel capacity in the finite-blocklength regime, which is  tighter than the famous closed-form asymptotic expression by  Polyanskiy-Poor-Verdú, given by
 \begin{align}\label{eq:Polyanskiy_rate}
 R(n, \epsilon) \approx C - \sqrt{\frac{V}{n}} Q^{-1}(\epsilon) + \frac{\log_2(n)}{2n},
 \end{align}
 where $C=\frac{1}{2}\log_2(1+\Omega)$  is the infinite block-length channel capacity, $n$ is the block/codeword length, $V=\frac{\Omega(\Omega+2)}{2(\Omega+1)^2} (\log_2 e)^2$,  $\Omega$ is the signal-to-noise-ratio (SNR),  $\epsilon$ is the error-rate, and  $R(n, \epsilon)$ is the data rate achieved with blocklength $n$ and error rate $\epsilon$. 
The expression in  \eqref{eq:Polyanskiy_rate} is obtained via normal approximation, see \cite{5452208}. Our novel   closed-form asymptotic expression of the bound on the AWGN channel capacity in the finite-blocklength regime is given by
 \begin{align}\label{eq:My_rate}
  R(n, \epsilon)  &\approx C-\frac{1}{2}\log_2\!\left(\frac{2\lambda}{ 1+\sqrt{1+\frac{4\lambda}{\Omega}}}\right) \nonumber\\
 &+\frac{1}{2\ln(2)}\left(1+\frac{1}{\Omega}+\frac{\lambda}{1+\Omega}- \sqrt{1+\frac{4\lambda}{\Omega}}\right)   ,
 \end{align}
 where
  \begin{align}\label{eq:My_rate_lambda}
 \lambda=1+\frac{1}{\Omega}+\frac{1}{\sqrt{n}} \sqrt{\frac{2(2+\Omega)}{\Omega} W\left(\frac{1}{2 \pi  \epsilon ^2}\right) },
 \end{align}
 where $W(\cdot)$ is the LambertW function. Note that  $ R(n, \epsilon)$, given by both in \eqref{eq:Polyanskiy_rate} and \eqref{eq:My_rate}, satisfies $\lim\limits_{n\to\infty} R(n, \epsilon)=C$, if $\epsilon>0$. As we will show using numerical examples, \eqref{eq:My_rate} is much tighter approximation to the converse bound than \eqref{eq:Polyanskiy_rate}.

 Finally, we present numerical examples where we show the application of the proposed method on bounding the tails of the   Gaussian, beta prime, and the non-central chi-squared RVs. Next, we show numerical examples that confirm the tightness of the derived lower and upper bounds on the converse bound of the AWGN channel capacity in the finite-blocklength regime.
We also show via numerical examples that \eqref{eq:My_rate} is much tighter approximation to the converse bound than \eqref{eq:Polyanskiy_rate}. 
 
 In general, this work provides compact, closed-form tail bounds that can be evaluated without numerical integration and that recover classical inequalities (e.g., Markov and Chernoff) as special cases. Beyond theory, these bounds quantify rare-event probabilities that underpin engineering tasks such as link error rates in communications, chance constraints in stochastic control, and reliability analysis.

 The paper is organized as follows. In Sec.~\ref{sec-2}, we provide some preliminary notations. In  Secs.~\ref{sec-3} and~\ref{sec-5}, we provide the bounding methods for the right and the left tails, respectively. In  Secs.~\ref{sec-4} and~\ref{sec-6}, we provide the iterative bounding methods for the right and the left tails, respectively. In  Sec.~\ref{sec-7}, we provide the convergence rates between the upper and lower bounds. In Sec.~\ref{sec-awgn}, we investigate the  AWGN capacity. In Sec.~\ref{sec-num}, we provide numerical results and in Sec.~\ref{sec-conc} we provide the conclusion. Finally, all the   proofs are provided in the Appendix. 

 \section{Preliminaries}\label{sec-2}

 Let $X$ be a continuous RV. Let $F (x)$ and $f (x)$ be the cumulative distribution function (CDF) and  the PDF of $X$, given by
 \begin{align}
 F (x)&={\rm Pr}\{X\leq  x\}\label{eq_1},\\
  f (x)&=\diff{ F (x)}{  x} ,\label{eq_2}
 \end{align}
 where  ${\rm Pr}\{A\}$ denotes the probability of an event $A$.

  The probabilities 
  \begin{align}
 {\rm Pr}\{X\leq  x\}=  F (x)\label{eq_n_1a} 
 \end{align}
 and
   \begin{align}
{\rm Pr}\{X\geq  x\}= 1-F (x) ,\label{eq_n_1b}
 \end{align}
 are known as the left tail and the  right tail probabilities of $X$, respectively.

In the following, we denote the first and the second derivatives of some function $y (x)$, by $y' (x)$ and $y'' (x)$, respectively. Thereby,  $f' (x)$ and $f'' (x)$ are the first and second derivatives of the PDF, $f (x)$, given by
  \begin{align}
  f' (x)&=\diff{  f (x)}{  x}  \label{eq_3}
 \end{align}
and
  \begin{align}
  f'' (x)&=\diff[2]{ f (x)}{  x} , \label{eq_4}
 \end{align} 
 respectively.
 
In the following, we assume that the PDF of the RV $X$, $f (x)$, has  support on  $[l,r]$, or on $(l,r]$, or on $[l,r)$, or on $(l,r)$, where $-\infty\leq  l< r \leq \infty$, which for simplicity we denote   as $[(l,r)]$. We assume that $f (x)$ is a continuous function of $x$ on the entire  support of $X$, i.e., that its  derivative    $f' (x)$ exists. Moreover, throughout this paper, when we write $\forall x>x_0$, we mean   $\forall x\in (x_0,r)$, and  when we write $\forall x<x_0$, we mean   $\forall x\in (l,x_0)$.

 \section{Bounds On The Right Tail}\label{sec-3}
 
In this section, we provide general upper and lower bounds on the right tail, $1-F (x)$, followed by a discussion about these bounds. We then provide some special cases. Finally, we connect the derived bounds to Markov's inequality and to Chernoff's bound.

\subsection{The General Bounds} 
 
 We start with the following useful lemma.
 
 \begin{lemma}\label{lema_1}
Let  $g (x)$  be any continuous, positive, and strictly decreasing function on a given  interval $I$, i.e., $g (x):\;$ $g (x)>0$ and $g' (x)< 0$, $\forall x\in I$. For such a $g (x)$, if
\begin{align}\label{seq_1_eq_1}
\frac{1-F (x)}{g (x)} 
\end{align}
is a deceasing function on the interval $I$, i.e., if the following holds
\begin{align}\label{seq_1_eq_2}
\diff{}{x}\left(\frac{1-F (x)}{g (x)} \right)\leq 0, \; \forall x\in I,
\end{align}
then the following upper holds
\begin{align}\label{seq_1_eq_3}
1-F (x)\leq - f (x)\frac{g (x)}{g' (x)} , \; \forall x\in I.
\end{align}
Otherwise,  if
\begin{align}\label{seq_1_eq_4}
\frac{1-F (x)}{g (x)} 
\end{align}
is an increasing function on the interval $I$, i.e., if the following holds
\begin{align}\label{seq_1_eq_5}
\diff{}{x}\left(\frac{1-F (x)}{g (x)} \right)\geq 0, \; \forall x\in I,
\end{align}
 then the following lower bound holds
\begin{align}\label{seq_1_eq_6}
1-F (x)\geq - f (x)\frac{g (x)}{g' (x)} , \; \forall x\in I.
\end{align}

 \end{lemma}

\begin{IEEEproof}
The proof is provided in Appendix~\ref{app_1}.
\end{IEEEproof}

Although the bounds in \eqref{seq_1_eq_3} and \eqref{seq_1_eq_6} seem simple, they are not practical since  determining whether condition \eqref{seq_1_eq_2} or condition \eqref{seq_1_eq_5} holds requires knowledge of $1-F (x)$, which by default we assume that is not available. Instead, we only know $f (x)$ and its derivative, $f' (x)$. This  practicality constraint is overcome by the following theorem, which provides  bounds similar to those in Lemma~\ref{lema_1}, but with corresponding conditions that depend only on $f (x)$ and $f' (x)$, and  not on $1-F (x)$.

\begin{theorem}\label{thm_1}
Let $P (x)$ be defined as
\begin{align}\label{seq_1_eq_7a}
P (x)=-   f (x)\frac{g (x)}{g' (x)},
\end{align}
where $g (x)$ is any continuous, positive, and strictly decreasing function  $\forall x>x_0$,  i.e., $g (x):\;$ $g (x)>0$, $g' (x)< 0$, $\forall x>x_0$. Moreover, let $g (x)$ be such that the following also holds
 \begin{align}\label{seq_1_eq_7b}
\lim_{x\to r} P (x)=  \lim_{x\to r}  - f (x)\frac{g (x)}{g' (x)}=0.
\end{align}

For any such function $g (x)$ as defined above, if
\begin{align}\label{seq_1_eq_8}
\frac{f (x)}{-g' (x)} 
\end{align}
is a deceasing function  $\forall x>x_0$, which is equivalent to the following condition being satisfied
\begin{align}\label{seq_1_eq_9}
P' (x)+f (x)\leq 0, \;  \forall x>x_0,
\end{align}
then  the following upper holds
\begin{align}\label{seq_1_eq_10}
1-F (x)\leq P (x) , \;  \forall x>x_0.
\end{align}
Otherwise,  if
\begin{align}\label{seq_1_eq_11}
\frac{f (x)}{-g' (x)} 
\end{align}
is an increasing function  $\forall x>x_0$,  which is equivalent to the following condition being satisfied
\begin{align}\label{seq_1_eq_12}
P' (x)+f (x)\geq 0, \;  \forall x>x_0,
\end{align}
 then  the following lower bound holds
\begin{align}\label{seq_1_eq_13}
1-F (x)\geq P (x) , \;  \forall x>x_0.
\end{align}
 \end{theorem}
 
\begin{IEEEproof}
The proof is provided in Appendix~\ref{app_2}.
\end{IEEEproof}

We now have a practical method to determine whether the bound in  \eqref{seq_1_eq_10} or the bound in \eqref{seq_1_eq_13} holds, simply by observing whether for a given $g (x)$, which satisfies the conditions defined in Theorem~\ref{thm_1}, condition \eqref{seq_1_eq_9} or condition \eqref{seq_1_eq_12} holds, respectively. Note that conditions   \eqref{seq_1_eq_9} and \eqref{seq_1_eq_12} depend only on $f (x)$, $f' (x)$, $g' (x)$, and $g'' (x)$, since 
\begin{align}
P' (x)&=\diff{}{x}\left(- f (x)\frac{g (x)}{g' (x)} \right)
\nonumber\\
&= 
- f' (x)\frac{g (x)}{g' (x)} - f (x) +  f (x)   \frac{g (x)g'' (x)}{\big(g' (x)\big)^2}.
\end{align}

In Theorem~\ref{thm_1}, note that we first need to provide a corresponding  function $g (x)$  and then check if it is a valid candidate for constructing an upper bound, a lower bound, or it is not a valid candidate.  There are many possible functions $g (x)$ that satisfy the conditions for $g (x)$ defined in  Theorem~\ref{thm_1}, and moreover satisfy either the upper bound condition in \eqref{seq_1_eq_9} or the lower bound condition in \eqref{seq_1_eq_12},  and thereby make the upper bound in \eqref{seq_1_eq_10} or  the lower bound in \eqref{seq_1_eq_13} to hold. 
But what is the optimal $g (x)$ for a given $f (x)$? It turns out that solving this problem, even for some special cases of $f (x)$, would require a standalone paper. Therefore, we leave the problem of finding the  optimal but practical $g (x)$, for a given $f (x)$, for future works. Note that the optimal but unpractical $g (x)$ always exists for a given $f (x)$, and is given by $g (x)=1-F (x)$, $\forall x\geq l$. If we plugin $g (x)=1-F (x)$, $\forall x\geq l$, into   Theorem~\ref{thm_1}, then it is easy to see that this $g (x)$ is a continuous, positive, and  decreasing\footnote{In the case when $g (x)=1-F (x)$, the strictly decreasing condition in Theorem~\ref{thm_1} can be replaced by  decreasing since for any point $l<x_0<r$ for which $\lim\limits_{x\to x_0}f (x)=0$, we have $\lim\limits_{x\to x_0} - f (x) \frac{g (x)}{g' (x)}= \lim\limits_{x\to x_0} - f (x) \frac{1-F (x)}{-f (x)}=1-F (x_0)$.} function that satisfies both the upper bound condition  in \eqref{seq_1_eq_9} and  the lower bound condition in \eqref{seq_1_eq_12}. Thereby, $P (x)$ constructed from $g (x)=1-F (x)$ is both an upper bound and a lower bound on $1-F (x)$, $\forall x\geq l$, which means that $P (x)=1-F (x)$.

For the problem of finding the  optimal but practical $g (x)$, for a given $f (x)$,  we    only provide the following intuitive observations.  A good upper bound $P (x)$ on  $1-F (x)$ is the one whose derivative $P' (x)$ is integrable in a closed-form expression,  $\forall x>x_0$, and $-P' (x)$ very tightly upper bounds $f (x)$, $\forall x>x_0$. The tighter $-P' (x)$   upper bounds $f (x)$, the tighter the upper bound $P (x)$ is on $1-F (x)$, $\forall x>x_0$. Similarly, a good lower bound $P (x)$ on  $1-F (x)$ is the one whose derivative $P' (x)$ is integrable in a  closed-form expression, $\forall x>x_0$,  and $-P' (x)$ very tightly lower bounds $f (x)$, $\forall x>x_0$. The tighter $-P' (x)$   lower bounds $f (x)$, the tighter the lower bound $P (x)$ is on $1-F (x)$, $\forall x>x_0$. In the limit, when $-P' (x)$   becomes equal to $f (x)$, $\forall x>x_0$, then $P (x)= 1-F (x)$, $\forall x>x_0$. However, we assume that, in this case, $P' (x)$  is not  integrable in a  closed-form expression, otherwise there won't be a need for bounding $1-F (x)$. The last claim  can be seen by solving the differential equation $P' (x)=-f (x)$, which results in $P (x)=1-F (x)$ and is obtained  by appropriately setting the  constant of the solution of this differential equation.

\begin{remark}\label{rem_1}
The gain in practicality provided by Theorem~\ref{thm_1} comes with a certain loss of generality as compared to Lemma~\ref{lema_1}.  For example, for a given $g (x)$ and for some distributions, the bounds  in Lemma~\ref{lema_1} hold $\forall x>x_L$, whereas Theorem~\ref{thm_1} shows that  the same bounds hold $\forall x>x_0$, where $x_L<x_0$. Thereby,  Theorem~\ref{thm_1} is ``blind'' to the fact that its bounds also hold in the interval $(x_L, x_0]$. We will encounter this situation later on when we establish a connection between the upper bound in Theorem~\ref{thm_1} and Markov's inequality.
 \end{remark}

 \subsection{Two Special Cases For $g (x)$} 
 
In general, there are many possible functions $g (x)$ that satisfy the conditions for $g (x)$ in  Theorem~\ref{thm_1}, and moreover satisfy either the upper bound condition in \eqref{seq_1_eq_9} or the lower bound condition in \eqref{seq_1_eq_12},  and thereby make the upper bound in \eqref{seq_1_eq_10} or  the lower bound in \eqref{seq_1_eq_13} to hold. In this subsection, we will  concentrate on  two such  functions for $g (x)$, which in many cases result in tight and/or simple upper bounds. These functions are given by
\begin{align} 
g (x)&=(x-l) f (x), \label{seq_1_eq_14a}\\
g (x)&=  f (x). \label{seq_1_eq_14b}  
\end{align}
 The upper bounds resulting from these two functions are provided in the following two corollaries.

\begin{cor}\label{cor_1}
Set $g (x)=(x-l) f (x)$. Let $g' (x)= f (x)+(x-l) f' (x)<0$, $\forall x>x_0$, and  
\begin{align}\label{seq_1_eq_15a}
\lim\limits_{x \to r}- f (x)\frac{(x-l) f (x)}{f (x)+(x-l) f' (x)} =0
\end{align}
 hold. If
\begin{align}\label{seq_1_eq_15}
-\frac{f (x)}{f (x)+(x-l) f' (x)} 
\end{align}
is a   deceasing function  $\forall x>x_0$, i.e., if   the following holds
\begin{align}\label{seq_1_eq_16}
\diff{ }{ x}\left(\frac{-f (x)}{f (x)+(x-l) f' (x)} \right)\leq 0, \;  \forall x>x_0,
\end{align}
then  the following upper holds
\begin{align}\label{seq_1_eq_17}
1-F (x)\leq - f (x)\frac{(x-l) f (x)}{f (x)+(x-l) f' (x)} , \;  \forall x>x_0.
\end{align}
 \end{cor}
 
 \begin{IEEEproof}
The proof is a direct result of  Theorem~\ref{thm_1}.
\end{IEEEproof}
 
 \begin{cor}\label{cor_2}
Set $g (x)= f (x)$. Let $g' (x)=  f' (x)<0$, $\forall x>x_0$, and  
\begin{align}\label{seq_1_eq_18a}
\lim\limits_{x\to r} - f (x)\frac{f (x)}{ f' (x)} =0
\end{align}
hold.  If
\begin{align}\label{seq_1_eq_18}
\frac{f (x)}{-f' (x)} 
\end{align}
is a deceasing function $\forall x>x_0$, i.e., if the following holds
\begin{align}\label{seq_1_eq_19}
\diff{ }{ x}\left(\frac{f (x)}{- f' (x)} \right)\leq 0, \;  \forall x>x_0,
\end{align}
then  the following upper holds
\begin{align}\label{seq_1_eq_20}
1-F (x)\leq - f (x)\frac{f (x)}{ f' (x)} , \;  \forall x>x_0.
\end{align}
 \end{cor}
 
  \begin{IEEEproof}
The proof is a direct result of   Theorem~\ref{thm_1}.
\end{IEEEproof}

The upper bounds in Corollaries~\ref{cor_1} and~\ref{cor_2} are simple  and yet tight for some of the most well known RVs, such as the Gaussian and  the chi-squared RVs, as will be shown in the numerical examples.

\subsection{Third Special Case For $g (x)$ and Connections to Markov's Inequality and Chernoff's Bound}
Another very special case for the function $g (x)$ is the following
\begin{align} 
g (x)&=\exp\left(-\int\limits_{x_0}^x \frac{f (z)}{h (z)} dz \right),\; \forall x>x_0, \label{seq_1_eq_14c} 
\end{align}
where   $h (x)$ will be defined in the following corollary.
For $g (x)$ given by \eqref{seq_1_eq_14c}, we have the following corollary.

\begin{cor}\label{cor_3}
For any continuous and positive function $h (x)$, $\forall x>x_0$, which   satisfies
\begin{align}\label{seq_1_eq_21}
\lim_{x\to r}   h (x) =0,
\end{align}
if 
\begin{align}\label{seq_1_eq_23}
h '(x)+f (x) \leq 0, \;  \forall x>x_0,
\end{align}
holds, then  the following upper holds
\begin{align}\label{seq_1_eq_24}
1-F (x)\leq h (x)  , \;  \forall x>x_0.
\end{align}
Otherwise,  if
\begin{align}\label{seq_1_eq_25}
h '(x)+f (x) \geq 0, \;  \forall x>x_0,
\end{align}
holds, then  the following lower bound holds
\begin{align}\label{seq_1_eq_27}
1-F (x)\geq h (x) , \;  \forall x>x_0.
\end{align}
 \end{cor}

\begin{IEEEproof}
Inserting $g (x)$ given by \eqref{seq_1_eq_14c} into Theorem~\ref{thm_1} and simplifying leads directly to this corollary. 
\end{IEEEproof}

We can now relate the bound in Corollary~\ref{cor_3} to Markov's inequality for RVs with unbounded support from the right. Specifically, let $X$ be a non-negative RV with  unbounded support from the right, i.e.,  $r\to\infty$ holds.  Then,  by setting $h (x)$  in Corollary~\ref{cor_3} as 
\begin{align}\label{seq_1_eq_28}
 h (x) =\frac{E\{X\}}{x}, \;  \forall x>0,
\end{align}
and by assuming  $0<E\{X\}<\infty$,
 we obtain the bound
 \begin{align}\label{seq_1_eq_29}
1-F (x)\leq \frac{E\{X\}}{x} , \;  \forall x>x_0,
\end{align}
where $x_0$ is the point for which \eqref{seq_1_eq_23}, i.e., the following begins to hold
\begin{align}\label{seq_1_eq_30}
-\frac{E\{X\}}{x^2}+f (x) \leq 0, \;  \forall x>x_0,
\end{align}
which is equivalent to
 \begin{align}\label{seq_1_eq_31}
 E\{X\} - x^2 f (x) \geq 0, \;  \forall x>x_0.
\end{align}
Although we know that Markov's inequality holds for all non-negative RVs and $\forall x\in [l,r]$,   Corollary~\ref{cor_3}    ``sees''  that  Markov's inequality holds for a) RVs with unbounded support from the right and b)  for $\forall x>x_0$, where $x_0$ is the point for which \eqref{seq_1_eq_31} begins to hold. This is because Corollary~\ref{cor_3} is ``blind'' to a)  functions $h (x)$ that do not satisfy \eqref{seq_1_eq_21}. Note that, for  $h (x)$ given by \eqref{seq_1_eq_28}, $h (x)$  satisfies \eqref{seq_1_eq_21} only if $E\{X\}<\infty$ and $r\to\infty$. Moreover, for b), Theorem~\ref{thm_1}, and thereby Corollary~\ref{cor_3}, are also ``blind' to   the interval $[l,x_0]$, as explained in Remark~\ref{rem_1}. However, Corollary~\ref{cor_3} also shows us a workaround for RVs with bounded support from the right, as explained in the following.

Let  $X$ be a non-negative RV with support on $([l,r])$. For such an RV, if we set $h (x)$  in Corollary~\ref{cor_3} as 
\begin{align}\label{seq_1_eq_28-1}
 h (x) =\frac{E\{X\}}{x} - \frac{E\{X\}}{r}, \;  \forall x>0,
\end{align}
and by assuming $0<E\{X\}<\infty$,
 we obtain the bound
 \begin{align}\label{seq_1_eq_29-1}
1-F (x)\leq \frac{E\{X\}}{x} - \frac{E\{X\}}{r}, \;  \forall x>x_0,
\end{align}
where $x_0$ is the point for which \eqref{seq_1_eq_23}, i.e., the following begins to hold
\begin{align}\label{seq_1_eq_30-1}
-\frac{E\{X\}}{x^2}+f (x) \leq 0, \;  \forall x>x_0,
\end{align}
which is equivalent to
 \begin{align}\label{seq_1_eq_31-1}
 E\{X\} - x^2 f (x) \geq 0, \;  \forall x>x_0.
\end{align}
Hence, by setting  $h (x)$ as in \eqref{seq_1_eq_28-1}, we obtain a bound similar to the one in Markov's inequality, but which is tighter than Markov's inequality for $0<r<\infty$ and for which we are certain that it holds for $x>x_0$.

Note also that Corollary~\ref{cor_3} provides a very general method for including the mean $E\{X\}$ into the tail bounds. Specifically,   by setting $h (x)$ in  Corollary~\ref{cor_3}  as 
\begin{align}\label{seq_1_cor_3-eq-1}
 h (x) =Q\big(x,  E\{X\} \big),
\end{align}
where   $Q\big(x,  E\{X\} \big)$ is any continuous function that satisfies $Q\big(x,  E\{X\} \big)>0$, $\forall x>x_0$, and $\lim\limits_{x\to r}Q\big(x,  E\{X\} \big)=0$, 
 we obtain the bound
 \begin{align}\label{seq_1_cor_3-eq-2}
1-F (x)\leq   Q\big(x,  E\{X\} \big) , \;  \forall x>x_0, 
\end{align}
when the following condition is met
 \begin{align}\label{seq_1_cor_3-eq-3}
   Q'\big(x,  E\{X\} \big)  +f (x) \leq 0, \;  \forall x>x_0.
\end{align}
Otherwise, we obtain the bound
 \begin{align}\label{seq_1_cor_3-eq-4}
1-F (x)\geq   Q\big(x,  E\{X\} \big) , \;  \forall x>x_0, 
\end{align}
when the following condition is met
 \begin{align}\label{seq_1_cor_3-eq-5}
   Q'\big(x,  E\{X\} \big)  +f (x) \geq 0, \;  \forall x>x_0.
\end{align}

It is straightforward to relate the bound in Corollary~\ref{cor_3} to other well known bounds such as the generalized Markov inequality, the Chebyshev's inequality, etc. 

Another case arises when we set an optimization parameter into the function $h(x)$ in Corollary~\ref{cor_3}, as specified in the following corollary.

\begin{cor}\label{cor_4}
For any continuous and positive   function $h \big(x,t(x)\big)$, $\forall x>x_0$, which also satisfies
\begin{align}\label{seq_1_eq_32}
\lim_{x\to r}   h \big(x,t(x)\big) =0,
\end{align}
if 
\begin{align}\label{seq_1_eq_33}
\frac{\partial h (x,t)}{\partial x} \Bigg|_{t=t(x)}+\frac{\partial h (x,t)}{\partial t}   \Bigg|_{t=t(x)} \frac{\partial t(x)}{\partial x} +f (x) \leq 0, \;  \forall x>x_0,
\end{align}
then  the following upper holds
\begin{align}\label{seq_1_eq_34}
1-F (x)\leq   h \big(x,t(x)\big)  , \;  \forall x>x_0,
\end{align}
Otherwise,  if
\begin{align}\label{seq_1_eq_35}
\frac{\partial h (x,t)}{\partial x} \Bigg|_{t=t(x)}+\frac{\partial h (x,t)}{\partial t}   \Bigg|_{t=t(x)} \frac{\partial t(x)}{\partial x} +f (x) \geq 0, \;  \forall x>x_0,
\end{align}
then  the following lower bound holds
\begin{align}\label{seq_1_eq_36}
1-F (x)\geq h \big(x,t(x)\big) , \;  \forall x>x_0.
\end{align}
 \end{cor}

 \begin{IEEEproof}
Replacing $h (x)$ with $h \big(x,t(x)\big)$ in Corollary~\ref{cor_3}, and then using the total derivative rule, leads directly to this corollary.
\end{IEEEproof}

Using Corollary~\ref{cor_4}, we can now construct optimization problems for tightening the tail bounds via the parameter $t$. Specifically, for the upper bound, the optimization problem would be
 \begin{align}\label{seq_1_nn1}
\min_t&\;\; h \big(x,t(x)\big),\;  \forall x>x_0 \nonumber\\
\textrm{s.t.} \;\;\;
&C_1: \; \eqref{seq_1_eq_33}
\nonumber\\
&C_2: \; h \big(x,t(x)\big)>0,\;  \forall x>x_0\nonumber\\
&C_3: \lim_{x\to r}   h \big(x,t(x)\big) =0.
\end{align}
For the lower bound, the optimization problem would be
 \begin{align}\label{seq_1_nn2}
\max_t&\; \;h \big(x,t(x)\big),\;  \forall x>x_0 \nonumber\\
\textrm{s.t.} \;\;\;
&C_1: \; \eqref{seq_1_eq_35}
\nonumber\\
&C_2: \; h \big(x,t(x)\big)>0,\;  \forall x>x_0\nonumber\\
&C_3: \lim_{x\to r}   h \big(x,t(x)\big) =0.
\end{align}

Moreover, we can now relate the bound in Corollary~\ref{cor_4} to the  Chernoff's bound for RVs with unbounded support from the right.   Specifically, by setting $h (x,t(x))$ in  Corollary~\ref{cor_4}  as 
\begin{align}\label{seq_1_eq_37}
 h (x,t(x)) =\min_t M (t)   e^{-t x},
\end{align}
where $M (t)$ is the moment generating function (MGF), we obtain the bound
 \begin{align}\label{seq_1_eq_38}
1-F (x)\leq \min_t M (t)   e^{-t x} , \;  \forall x>x_0,
\end{align}
where $x_0$ is the point for which \eqref{seq_1_eq_33}, i.e., the following begins to hold
\begin{align}\label{seq_1_eq_39}
&-t M (t)     e^{-t x}   \Bigg|_{t=t^*(x)} \hspace{-2mm}+ \left(M' (t)   e^{-t x} - x M (t)   e^{-t x} \right)  \Bigg|_{t=t^*(x)}  {t^*}' (x)\nonumber\\
& +f (x) \leq 0, \;  \forall x>x_0,   
\end{align}
where $t^*(x)$ is the solution to $\min\limits_t M (t)   e^{-t x}$.

For RVs with  bounded support from the right, similar to the Markov's inequality type of bound explained above, we can set  $h (x,t(x))$ in  Corollary~\ref{cor_4}  as 
\begin{align}\label{seq_1_eq_39-1}
 h (x,t(x)) =\min_t M (t)  \big( e^{-t x} -  e^{-t r}\big).
\end{align}

Finally, the inclusion of the MGF into the tail bounds can be made in a much more general manner by setting $h (x,t(x))$ in  Corollary~\ref{cor_4}  as 
\begin{align}\label{seq_1_eq_39-2}
 h (x,t(x)) =Q\big(x, t, M (t) \big),
\end{align}
where   $Q\big(x, t, M (t) \big)$ is any continuous function that satisfies $Q\big(x, t, M (t) \big)>0$, $\forall x>x_0$, and  
$$\lim\limits_{x\to r}Q\big(x, t, M (t) \big)=0.$$ Thereby,  inserting \eqref{seq_1_eq_39-2} into the optimization problems in \eqref{seq_1_nn1} and \eqref{seq_1_nn2}, would result in upper and lower bounds as functions of the MGF, respectively.

We have yet to provide functions $g (x)$ that are good candidates for the corresponding lower bound in Theorem~\ref{thm_1}. Such $g (x)$ will be arrived  at by an iterative method, which is the subject of the following section.

\section{The Iterative Method For The Right Tail}\label{sec-4}

In this section, we provide an iterative method for obtaining ever tighter upper and lower bounds on $1-F (x)$, under certain conditions. Before we provide the iterative method, we introduce several lemmas which will be useful  for arriving at the iterative method. Moreover, in this section, when we say that some function is a valid upper or a lower bound on $1-F (x)$ as per Theorem~\ref{thm_1}, we mean that this bound is obtained using Theorem~\ref{thm_1} and thereby satisfies all of the conditions laid out in Theorem~\ref{thm_1}.

We start with the following lemma.
 \begin{lemma}\label{lema_2}
If $g (x)$, with properties defined in Theorem~\ref{thm_1}, also satisfies the  limit  
 \begin{align}\label{seq_1_eq_13-n1}
 \lim_{x\to r} g (x)=0,
 \end{align}
 then $P (x)$, given by \eqref{seq_1_eq_7a},  satisfies \eqref{seq_1_eq_7b}, i.e.,  $P (x)$ satisfies the limit  
 \begin{align}\label{seq_1_eq_13-n2}
\lim_{x\to r} P (x)=  \lim_{x\to r}  -  f (x)\frac{g (x)}{g' (x)}=0.
\end{align}
 
 \end{lemma}
\begin{IEEEproof}
The proof is provided in Appendix~\ref{app_2-3}.
\end{IEEEproof}

We note that a function $g (x)$ does not need to satisfy   \eqref{seq_1_eq_13-n1} in order for \eqref{seq_1_eq_7b} (i.e.,  \eqref{seq_1_eq_13-n2}) to hold. In other words, there are functions $g (x)$ for which \eqref{seq_1_eq_13-n1} does not hold and yet \eqref{seq_1_eq_7b} (i.e.,  \eqref{seq_1_eq_13-n2}) holds. However, what Lemma~\ref{lema_2} shows  us is that if $g (x)$ is such that \eqref{seq_1_eq_13-n1} holds, then we have certainty that \eqref{seq_1_eq_7b} (i.e.,  \eqref{seq_1_eq_13-n2})   holds. We will find   Lemma~\ref{lema_2} useful later on.

We now start providing the basic building elements of the iterative method. 

Let us define $P_0(x)$ as
\begin{align}\label{seq_2_eq_1}
P_0(x)=-   f (x)\frac{g (x)}{g' (x)}.
\end{align}
Note that $P_0(x)$ in \eqref{seq_2_eq_1} is identical to $P (x)$ given by \eqref{seq_1_eq_7a}.
Let us assume that $g (x)$ in \eqref{seq_2_eq_1} satisfies the conditions defined in Theorem~\ref{thm_1} in order for \eqref{seq_2_eq_1} to be an upper on $1-F (x)$, $\forall x>x_0$, or for  \eqref{seq_2_eq_1} to be a lower on $1-F (x)$, $\forall x>x_0$.
 Next, let us define the function $P_i(x)$, for $i=0,1,2,...$ as
 \begin{align}\label{seq_2_eq_2}
P_{i+1}(x)=-   f (x)\frac{P_i(x)}{P'_i(x)}.
\end{align}
Note that $P_{i+1}(x)$ is obtained in an itterative manner starting from the seed $P_{0}(x)$ given by \eqref{seq_2_eq_1}.

For the limit of $P_{i+1}(x)$ as $x\to r$, we have the following lemma.

 \begin{lemma}\label{lema_2a}
If $P_0(x)$ is a  valid upper bound or a valid lower bound  on $1-F (x)$, $\forall x>x_0$, as per Theorem~\ref{thm_1}, then the following limit holds for any $i=0,1,2,...,$
 \begin{align}\label{seq_1_eq_13-n2-1}
\lim\limits_{x\to r} P_{i+1}(x)=0.
\end{align}
 \end{lemma}
 
 \begin{IEEEproof}
The proof is provided in Appendix~\ref{app_2-3a}.
\end{IEEEproof}

Next, we have the following useful lemma for $P_{i+1}(x)$.

\begin{lemma}\label{lema_3}
If $P_{i}(x)$ is an upper on $1-F (x)$, $\forall x>x_0$, as per Theorem~\ref{thm_1}, which also satisfies $P'_{i}(x)<0$,  $\forall x>\hat x_0$, then the following holds
 \begin{align}\label{seq_2_eq_3}
P_{i+1}(x)\leq P_{i}(x), \; \forall x>x_i \textrm{ and any } i=0,1,2,... ,
\end{align}
where $x_i=\max(x_0,\hat x_0)$.

Otherwise, if $P_{i}(x)$ is a lower bound on $1-F (x)$, $\forall x>x_0$, as per Theorem~\ref{thm_1}, which also satisfies $P'_{i}(x)<0$,  $\forall x>\hat x_0$,    then the following holds
 \begin{align}\label{seq_2_eq_4}
P_{i+1}(x)\geq P_{i}(x),\; \forall x>x_i \textrm{ and any } i=0,1,2,...,
\end{align}
where $x_i=\max(x_0,\hat x_0)$.
\end{lemma}

 \begin{IEEEproof}
The proof is provided in Appendix~\ref{app_3}.
\end{IEEEproof}

Lemma~\ref{lema_3} is useful since it tells us that the next iteration  $P_{i+1}(x)$ of an upper bound $P_i(x)$ is always smaller than the preceding iteration  $P_i(x)$. Thereby, if the next iteration,  $P_{i+1}(x)$, itself is also an upper bound, then  $P_{i+1}(x)$ will be a tighter upper bound than its preceding iteration,  $P_i(x)$. Similarly, Lemma~\ref{lema_3}  tells us that the next iteration,  $P_{i+1}(x)$, of a lower bound $P_i(x)$ is always larger than the preceding iteration, $P_i(x)$. Thereby, if the new iteration, $P_{i+1}(x)$, itself is also a lower bound, then $P_{i+1}(x)$ will be a tighter lower bound than its preceding iteration, $P_{i}(x)$. We can imagine that problems may arise if the next iteration of an upper bound (lower bound) becomes a lower bound (upper bound). But when that happens, we have methods to check if that will result in a tighter bound than the one from the previous iteration, as explained in the following.

We now introduce  an auxiliary lower bound which we can use to measure the tightness of a lower bound obtained iteratively from a preceding upper bound. Specifically,
 note that from a given upper bound, $P_i(x)$, we can always create an auxiliary  lower bound, denoted by $P_{L,i}(x)$, which is obtained by reflecting $P_i(x)$ with respect to $1-F (x)$, as
 \begin{align}\label{seq_2_eq_4a}
P_{L,i}(x)&=  1-F (x)-\Big(P_i(i)-\big(1-F (x)\big)\Big) 
\nonumber\\
&=  2(1-F (x)) - P_i(x), \; \forall x>x_i,
\end{align}
when $P_i(x)$ is a valid upper  bound on $1-F (x)$  as per Theorem~\ref{thm_1}.

Similarly, we introduce  an auxiliary upper bound which we can use to measure the tightness of an upper bound obtained iteratively from a preceding lower bound. Again, note that from a given lower bound, $P_i(x)$, we can always create an auxiliary  upper bound, denoted by $P_{U,i}(x)$, which is obtained by reflecting $P_i(x)$ with respect to $1-F (x)$, as
 \begin{align}\label{seq_2_eq_4aa}
P_{U,i}(x) &=  1-F (x)+\Big(\big(1-F (x)\big)-P_i(i)\Big) 
\nonumber\\
&=  2(1-F (x)) - P_i(x), \; \forall x>x_i,
\end{align}
when $P_i(x)$ is a valid lower  bound on $1-F (x)$  as per Theorem~\ref{thm_1}.

Using $P_{L,i}(x)$ and $P_{U,i}(x)$, we can state the following lemma.

\begin{lemma}\label{lema_4}
Let $P_{i}(x)$ be an upper on $1-F (x)$, $\forall x>x_i$, and let
 \begin{align}\label{seq_2_eq_5}
P_{i+1}(x)=-   f (x)\frac{P_i(x)}{P'_i(x)}, \; \forall x>x_{i+1},
\end{align}
be a lower bound on $1-F (x)$, $\forall x>x_{i+1}$. In that case, if the condition
\begin{align}\label{seq_2_eq_7}
P'_{i+1}(x)+P'_{i}(x)+2 f (x)\leq 0, \; \forall x>\hat x_{i+1},
\end{align}
holds, then 
$P_{i+1}(x)$ is a tighter bound on $1-F (x)$ than $P_{i}(x)$, in the sense that
 \begin{align}\label{seq_2_eq_8}
P_{L,i}(x)\leq  P_{i+1}(x), \; \forall \hat x>x_{i+1},
\end{align}
holds, where $P_{L,i}(x)$ is given by \eqref{seq_2_eq_4a}.

On the other hand, let $P_{i}(x)$ be a lower bound on $1-F (x)$, $\forall x>x_i$, and let
 \begin{align}\label{seq_2_eq_9}
P_{i+1}(x)=-   f (x)\frac{P_i(x)}{P'_i(x)}, \; \forall x>x_{i+1},
\end{align}
be an upper bound on $1-F (x)$, $\forall x>x_{i+1}$. In that case, if the condition
 \begin{align}\label{seq_2_eq_11}
P'_{i+1}(x)+P'_{i}(x)+2 f (x)\geq 0, \; \forall \hat x>x_{i+1},
\end{align}
holds, then $P_{i+1}(x)$ is a tighter bound on $1-F (x)$ than $P_{i}(x)$, in the sense that
 \begin{align}\label{seq_2_eq_10}
P_{U,i}(x)\geq  P_{i+1}(x), \; \forall \hat x>x_{i+1},
\end{align}
 where $P_{U,i}(x)$ is given by \eqref{seq_2_eq_4aa}.

\end{lemma}

 \begin{IEEEproof}
The proof is provided in Appendix~\ref{app_4}.
\end{IEEEproof}

We now have all of the necessary elements to provide an iterative algorithm that can lead to ever tighter upper and lower bounds obtained in an iterative manner; bounds which not necessary all start to hold from the same $x=x_0$.  Specifically, if $P_{i}(x)$ is an upper/lower bound on $1-F (x)$, $\forall x>x_i$, then it may happen that $P_{i+1}(x)$ is a lower/upper bound on $1-F (x)$, $\forall x>x_{i+1}$, where $x_{i+1}>x_i$. A note that for the Gaussian RV, the author has not been able to find a promising function $g (x)$ from which the iterative method can be started such that each next iteration is a valid upper/lower bound $\forall x>x_0$. Instead, in the author's experiments, each next iterative bound starts to hold for ever larger $x$.

\begin{algorithm}
\caption{The Iterative Method For The Right Tail}\label{alg:1}
 
\SetKwData{Left}{left}
\SetKwInOut{Set}{Set}
	\SetKwInOut{Ensure}{Ensure}
	\SetKwInOut{Input}{Input}
	\SetKwInOut{Output}{Output}
	\SetKwInOut{Return}{Return}
	\SetKwData{And}{and}
	\Input{$f (x)$, $g (x)$, $x_0$}
	\Output{$P_L(x)$, $P_U(x)$}
	\BlankLine
	\Set{$P_0(x)=-   f (x)\dfrac{g (x)}{g' (x)}$}
	\BlankLine
	\Set{$C=0$, $i=0$, $P_L(x)=NaN$, $P_U(x)=NaN$}
	\While{$C=0$  }{
	$P_{i+1}(x)=-   f (x)\dfrac{P_{i}(x)}{P'_{i}(x)}$\\
	\uIf{$P_{i}(x)> 0, \; \forall x>x_0$, \textbf{ and } $P'_{i}(x)< 0, \; \forall x>x_0$}
	{
	\uIf{$P'_{i}(x)+f (x)\leq 0, \; \forall x>x_0$,}
	{
	\uIf{$P'_{i+1}(x)+f (x)\leq 0, \; \forall x>x_0$,}
	{	\BlankLine
	$P_U(x)=P_{i+1}(x)$
		\BlankLine}
	\uElseIf{$P'_{i+1}(x)+f (x)\geq 0, \; \forall x>x_0$, \textbf{ and } $P'_{i+1}(x)+P'_{i}(x)+2f (x)\leq 0, \; \forall x>x_0$,}
	{	\BlankLine
	$P_L(x)=P_{i+1}(x)$
		\BlankLine}
	\Else{$C=1$}
	}
	\Else{
	\uIf{$P'_{i+1}(x)+f (x)\leq 0, \; \forall x>x_0$, \textbf{ and } $P'_{i+1}(x)+P'_{i}(x)+2f (x)\geq 0, \; \forall x>x_0$,}
	{$P_U(x)=P_{i+1}(x)$}
	\uElseIf{$P'_{i+1}(x)+f (x)\geq 0, \; \forall x>x_0$,}
	{$P_L(x)=P_{i+1}(x)$}
	\Else{$C=1$}
	}
	}
	\Else{$C=1$}
		$i=i+1$
	}
	\Return{$P_L(x)$, $P_R(x)$}
\end{algorithm}

The iterative algorithm is given in Algorithm~\ref{alg:1}, and works as follows.  The algorithm takes as inputs the PDF, $f (x)$, the function $g (x)$, and a desired point $x_0$ from which we want these iteratively obtained bounds to hold. Note, the function $g (x)$ must be such that $P_0(x)$ is a valid upper or a lower bound $\forall x>x_0$, as per Theorem~\ref{thm_1}. The \textit{while loop}  performs the following computations in an iterative manner, unless in the process of iteration  $C$ changes value from $C=0$ to $C=1$.
The   \textit{outer if  condition}, checks if the function $g (x)$ in iteration $i$, which in this case is $g (x)=P_i(x)$, is a valid function according to Theorem~\ref{thm_1}. If true, then the algorithm continues to the \textit{middle if  condition}. If false, $C$ is set to $C=1$ and the  \textit{while loop} stops.
The   \textit{middle if  condition} in the algorithm checks  if the previous iteration $P_i(x)$ is a valid upper bound. If true, then the \textit{inner if  condition}   checks if the next iteration $P_{i+1}(x)$ is also a valid  upper bound. If true, then $P_{i+1}(x)$ is a tighter upper bound than $P_i(x)$, according to Lemma~\ref{lema_3},  and therefore   $P_{i+1}(x)$ is stored into $P_U(x)$. Otherwise, if the \textit{inner if  condition} is false,  it is checked whether the next iteration $P_{i+1}(x)$ is  a valid  lower bound and if $P_{i+1}(x)$ is a tighter bound  than $P_{i}(x)$, as per Lemma~\ref{lema_4}. If true, then   $P_{i+1}(x)$ is stored into $P_L(x)$. If false, $C$ is set to $C=1$ and the \textit{while loop} stops.
 On the other hand, if the \textit{middle if  condition} is false, then the previous iteration $P_i(x)$ must be a valid lower bound. Therefore,  the \textit{inner if  condition} checks if the next iteration $P_{i+1}(x)$ is   a valid  upper bound and if this bound is tighter than $P_{i}(x)$, as per Lemma~\ref{lema_4}. If true, then $P_{i+1}(x)$ is stored into $P_U(x)$. If false,  then the next iteration $P_{i+1}(x)$ is checked whether it is valid  lower bound. If true, then $P_{i+1}(x)$ is a tighter lower bound than $P_i(x)$, according to Lemma~\ref{lema_3},  and therefore   $P_{i+1}(x)$  is stored   into $P_L(x)$. If false,    $C$ is set to $C=1$ and the \textit{while loop} stops. If $C$ has not changed value during one cycle in the \textit{while loop}, then $C$  keeps the value  $C=0$ and therefore the \textit{while loop} performs another iteration. 
 Finally, when the \textit{while loop} stops, the algorithm returns $P_L(x)$ and $P_U(x)$. 
 
 For a better understanding of the iterative method, we  provide an example for the Gaussian RV in the following.

\begin{example}
 Let $f (x)$ be the PDF of the zero-mean unit-variance Gaussian RV. Let us choose $g (x)$ as $g (x)=f (x)$. Note that $g (x)$ is continuous,  $g(x)>0$, $\forall x$, $g' (x)<0$, $\forall x>0$, and $\lim\limits_{x\to \infty}  - f (x)\frac{g (x)}{g' (x)}=0$. Hence, $g (x)$ satisfies all of the conditions laid out in Theorem~\ref{thm_1}. Using $g (x)$, we construct  $P_0(x)$ as in \eqref{seq_2_eq_1}, and thereby we obtain $P_0(x)$ as
\begin{align}\label{seq_2_eq_12}
P_0(x) =- f (x)\frac{f (x)}{f '(x)}=\frac{1}{x}\frac{e^{-\frac{x^2}{2}}}{ \sqrt{2 \pi } }.
\end{align}
It is easy to verify that \eqref{seq_2_eq_12} satisfies \eqref{seq_1_eq_9}, $\forall x$, and thereby \eqref{seq_2_eq_12} is a valid upper bound on $1-F (x)$,   $\forall x>0$ (the condition $\forall x>0$ comes from the fact that $g' (x)>0$, $\forall x>0$).

Next, from  $P_0(x)$ in \eqref{seq_2_eq_12}, we construct  $P_1(x)$ using \eqref{seq_2_eq_2}, and thereby we obtain $P_1(x)$ as
\begin{align}\label{seq_2_eq_13}
P_1(x) &=- f (x)\frac{P_0(x)}{P_0'(x)}= \frac{f^2 (x) f '(x)}{f (x) f ''(x)-2 \big(f '(x)\big)^2} \nonumber\\
&= \frac{x}{1+x^2}\frac{e^{-\frac{x^2}{2}} x}{\sqrt{2 \pi }}.
\end{align}
It is easy to verify that \eqref{seq_2_eq_13} satisfies \eqref{seq_1_eq_12}, $\forall x$, and thereby \eqref{seq_2_eq_13} is a valid lower bound on $1-F (x)$, $\forall x>0$ (the condition $\forall x>0$ now comes from the fact that  $P_0(x)>0$, $\forall x>0$). If we now check condition \eqref{seq_2_eq_7}, it is easy to verify that the lower bound in \eqref{seq_2_eq_13} is a tighter\footnote{Tighter in the sense of Lemma~\ref{lema_4}.\label{note1}} bound on $1-F (x)$ than  its preceding iteration, the upper bound in \eqref{seq_2_eq_12}, $\forall x>0$. 

If we do another iteration, now from  $P_1(x)$ in \eqref{seq_2_eq_13}, by constructing  $P_2(x)$ using \eqref{seq_2_eq_2}, we obtain $P_2(x)$ as
\begin{align}\label{seq_2_eq_14}
P_2(x) =- f (x)\frac{P_1(x)}{P_1'(x)}= \frac{ \left(x^3+x\right)}{\left(x^4+2 x^2-1\right)} \frac{e^{-\frac{x^2}{2}} }{\sqrt{2 \pi } }.
\end{align}

It is easy to verify that \eqref{seq_2_eq_14} satisfies \eqref{seq_1_eq_9}, $\forall x$, and thereby \eqref{seq_2_eq_14} is a valid upper bound on $1-F (x)$, $\forall x>\sqrt{\sqrt{2}-1}$ (the condition $\forall x>\sqrt{\sqrt{2}-1}$ now comes from the fact that  $P'_1(x)<0$, $\forall x>\sqrt{\sqrt{2}-1}$). If we now check condition \eqref{seq_2_eq_11}, it is easy to verify that the upper bound in \eqref{seq_2_eq_14} becomes a tighter\footref{note1} bound on $1-F (x)$ than  its preceding iteration, the lower bound in \eqref{seq_2_eq_13}, $\forall x> 1.71298$.

If we  continue further with the iterative method, for each iteration $i+1$, we will obtain a tighter\footref{note1} upper bound, $P_{i+1}(x)$ (if its predecessor $P_i(x)$ was a lower bound) or we will obtain a tighter\footref{note1} lower bound, $P_{i+1}(x)$  (if its predecessor $P_i(x)$ was an upper bound), but  in  each second iteration, these bounds will hold for ever larger $x$. We will also see this property via numerical examples in Sec.~\ref{sec-num}. This ends this example.

\end{example}

%%%%%%%%%%%%%%%%%%%%%%%%%%%%%%%%%%%%%%%%%%%%%%%%%%%%%%%%%%%%%%%%%%%%%%%
%%%%%%%%%%%%%%%%%%%%%%%%%%%%%%%%%%%%%%%%%%%%%%%%%%%%%%%%%%%%%%%%%%%%%%%

 \section{Bounds On The Left Tail}\label{sec-5}
 
We now provide a mirror like results for the left tail.

\subsection{The General Bounds} 
 
 We start directly with the main theorem.

\begin{theorem}\label{thm_2}
Let $P (x)$ be defined as
\begin{align}\label{seq_5_eq_7a}
P (x)=   f (x)\frac{g (x)}{g' (x)},
\end{align}
where $g (x)$ is any continuous, positive, and strictly increasing function  $\forall x<x_0$,  i.e., $g (x):\;$ $g (x)>0$, $g' (x)> 0$, $\forall x<x_0$. Moreover, let $g (x)$ be such that the following also holds
 \begin{align}\label{seq_5_eq_7b}
\lim_{x\to l} P (x)=  \lim_{x\to l}   f (x)\frac{g (x)}{g' (x)}=0.
\end{align}

For any such function $g (x)$ as defined above, if
\begin{align}\label{seq_5_eq_8}
\frac{f (x)}{g' (x)} 
\end{align}
is an increasing function  $\forall x<x_0$, which is equivalent to the following condition being satisfied
\begin{align}\label{seq_5_eq_9}
P' (x)-f (x)\geq 0, \;  \forall x<x_0,
\end{align}
then  the following upper holds
\begin{align}\label{seq_5_eq_10}
 F (x)\leq P (x) , \;  \forall x<x_0.
\end{align}
Otherwise,  if
\begin{align}\label{seq_5_eq_11}
\frac{f (x)}{g' (x)} 
\end{align}
is a  decreasing function  $\forall x<x_0$,  which is equivalent to the following condition being satisfied
\begin{align}\label{seq_5_eq_12}
P' (x)-f (x)\leq 0, \;  \forall x<x_0,
\end{align}
 then  the following lower bound holds
\begin{align}\label{seq_5_eq_13}
F (x)\geq P (x) , \;  \forall x<x_0.
\end{align}
 \end{theorem}
 
\begin{IEEEproof}
The proof is provided in Appendix~\ref{app_5}.
\end{IEEEproof}

For the left tail, we now also have a practical method to determine whether the bound in  \eqref{seq_5_eq_10} or the bound in \eqref{seq_5_eq_13} holds, simply by observing whether for a given $g (x)$, which satisfies the conditions defined in Theorem~\ref{thm_2}, condition \eqref{seq_5_eq_9} or condition \eqref{seq_5_eq_12} holds, respectively.

 \subsection{A Special Case  For \lowercase{$g (x)$}} 
 
There are many possible functions $g (x)$ that satisfy the conditions for $g (x)$ defined in  Theorem~\ref{thm_2}, and moreover satisfy either the upper bound condition in \eqref{seq_5_eq_9} or the lower bound condition in \eqref{seq_5_eq_12},  and thereby make the upper bound in \eqref{seq_5_eq_10} or  the lower bound in \eqref{seq_5_eq_13} to hold. In this subsection, we will  concentrate on  one such  function, $g (x)$, which in many cases result in a tight and/or simple upper bound. This function  is given by
\begin{align} 
g (x)&=(x-l) f (x), \label{seq_5_eq_14a}
\end{align}
 The upper bound  resulting from this   function  is provided in the following  corollary.

\begin{cor}\label{cor_5}
Set $g (x)=(x-l) f (x)$. Let $g' (x)= f (x)+(x-l) f' (x)>0$, $\forall x<x_0$, and  
\begin{align}\label{seq_5_eq_15a}
\lim\limits_{x \to l} f (x)\frac{(x-l) f (x)}{f (x)+(x-l) f' (x)} =0
\end{align}
 hold. If
\begin{align}\label{seq_5_eq_15}
 \frac{f (x)}{f (x)+(x-l) f' (x)} 
\end{align}
is an   increasing function  $\forall x<x_0$, i.e., if   the following holds
\begin{align}\label{seq_5_eq_16}
\diff{ }{ x}\left(\frac{ f (x)}{f (x)+(x-l) f' (x)} \right)\geq 0, \;  \forall x<x_0,
\end{align}
then  the following upper holds
\begin{align}\label{seq_5_eq_17}
 F (x)\leq  f (x)\frac{(x-l) f (x)}{f (x)+(x-l) f' (x)} , \;  \forall x<x_0.
\end{align}
 \end{cor}
 
 \begin{IEEEproof}
The proof is a direct result of  Theorem~\ref{thm_2}.
\end{IEEEproof}

The upper bound  in Corollary~\ref{cor_5} is simple  and yet tight for some of the most well known RVs, such as the  the chi-squared RV, as will be shown in the numerical examples.

\subsection{Another Special Case For \lowercase{$g (x)$} and Constructing Left Tail Bounds Similar to Markov's Inequality and Chernoff's Bound}
Another very special case for the function $g (x)$ is the following
\begin{align} 
g (x)&=\exp\left(\int\limits_{l}^x \frac{f (z)}{h (z)} dz \right), \; \forall x<x_0, \label{seq_5_eq_14c} 
\end{align}
where  $h (x)$ will be defined in the following corollary.
For $g (x)$ given by \eqref{seq_5_eq_14c}, we have the following corollary.

\begin{cor}\label{cor_6}
For any continuous and positive   function $h (x)$, $\forall x<x_0$, which also satisfies
\begin{align}\label{seq_5_eq_21}
\lim_{x\to l}   h (x) =0,
\end{align}
if 
\begin{align}\label{seq_5_eq_23}
h '(x)-f (x) \geq 0, \;  \forall x<x_0,
\end{align}
then  the following upper holds
\begin{align}\label{seq_5_eq_24}
 F (x)\leq h (x)  , \;  \forall x<x_0.
\end{align}
Otherwise,  if
\begin{align}\label{seq_5_eq_25}
h '(x)-f (x) \leq 0, \;  \forall x<x_0,
\end{align}
then  the following lower bound holds
\begin{align}\label{seq_5_eq_27}
 F (x)\geq h (x) , \;  \forall x<x_0.
\end{align}
 \end{cor}

\begin{IEEEproof}
Inserting $g (x)$ given by \eqref{seq_5_eq_14c} into Theorem~\ref{thm_2} and simplifying leads directly to this corollary. 
\end{IEEEproof}

We know that Markov's inequality  holds for the right tail only. However, 
using Corollary~\ref{cor_6}, we can now create a type of left tail bounds similar to Markov's inequality, in the sense that the mean will be included into the tail bound. Specifically, let $X$ be an RV with    support  on $([l,r])$.    Then,  by setting $h (x)$  in Corollary~\ref{cor_6} as 
\begin{align}\label{seq_5_eq_28}
 h (x) = Q(x, E\{X\}) , \;  \forall x\leq x_0,
\end{align}
where $Q(x, E\{X\})$ is any continuous function that satisfies $Q(x, E\{X\})>0$ and  $\lim\limits_{x\to l}Q(x, E\{X\})=0$, 
 we obtain the bound
 \begin{align}\label{seq_5_eq_29}
F (x)\leq  Q(x, E\{X\}) , \;  \forall x<x_0,
\end{align}
when \eqref{seq_5_eq_23}, i.e., the following condition  holds
\begin{align}\label{seq_5_eq_30}
Q'(x, E\{X\}) -f (x) \geq 0, \;  \forall x<x_0.
\end{align}
On the other hand, the following lower bound holds
 \begin{align}\label{seq_5_eq_30-1}
F (x)\geq  Q(x, E\{X\}) , \;  \forall x<x_0,
\end{align}
when \eqref{seq_5_eq_25}, i.e., the following condition holds
\begin{align}\label{seq_5_eq_30-2}
Q'(x, E\{X\}) -f (x) \leq 0, \;  \forall x<x_0.
\end{align}

Another case arises when we set an optimization parameter into the function $h(x)$ in Corollary~\ref{cor_6}, as specified in the following corollary.

\begin{cor}\label{cor_7}
For any continuous and positive  function $h \big(x,t(x)\big)$, $\forall x<x_0$,  which also satisfies
\begin{align}\label{seq_5_eq_32}
\lim_{x\to l}   h \big(x,t(x)\big) =0,
\end{align}
if 
\begin{align}\label{seq_5_eq_33}
\frac{\partial h (x,t)}{\partial x} \Bigg|_{t=t(x)} \hspace{-2mm}+\frac{\partial h (x,t)}{\partial t}   \Bigg|_{t=t(x)} \frac{\partial t(x)}{\partial x} -f (x) \geq 0, \;  \forall x<x_0,
\end{align}
then  the following upper holds
\begin{align}\label{seq_5_eq_34}
 F (x)\leq   h \big(x,t(x)\big)  , \;  \forall x<x_0,
\end{align}
Otherwise,  if
\begin{align}\label{seq_5_eq_35}
\frac{\partial h (x,t)}{\partial x} \Bigg|_{t=t(x)}\hspace{-2mm}+\frac{\partial h (x,t)}{\partial t}   \Bigg|_{t=t(x)} \frac{\partial t(x)}{\partial x} -f (x) \leq 0, \;  \forall x<x_0,
\end{align}
then  the following lower bound holds
\begin{align}\label{seq_5_eq_36}
 F (x)\geq h \big(x,t(x)\big) , \;  \forall x<x_0.
\end{align}
 \end{cor}

 \begin{IEEEproof}
Replacing $h (x)$ with $h \big(x,t(x)\big)$ in Corollary~\ref{cor_6}, and then using the total derivative rule, leads directly to this corollary.
\end{IEEEproof}

Using Corollary~\ref{cor_7}, we can now construct optimization problems for tightening the bounds. Specifically, for the upper bound, the optimization problem would be
 \begin{align}\label{seq_5_nn1}
\min_t&\;\; h \big(x,t(x)\big),\;  \forall x<x_0 \nonumber\\
\textrm{s.t.} \;\;\;
&C_1: \; \eqref{seq_5_eq_33}
\nonumber\\
&C_2: \; h \big(x,t(x)\big)>0,\;  \forall x<x_0\nonumber\\
&C_3: \; \lim_{x\to l}   h \big(x,t(x)\big) =0.
\end{align}
For the lower bound, the optimization problem would be
 \begin{align}\label{seq_5_nn2}
\max_t&\; \;h \big(x,t(x)\big),\;  \forall x<x_0 \nonumber\\
\textrm{s.t.} \;\;\;
&C_1: \; \eqref{seq_5_eq_35}
\nonumber\\
&C_2: \; h \big(x,t(x)\big)>0,\;  \forall x<x_0\nonumber\\
&C_3: \lim_{x\to l}   h \big(x,t(x)\big) =0.
\end{align}

Using Corollary~\ref{cor_7}, we can now construct  left tail bounds similar to  Chernoff's bound, in the sense that the MGF will be included into the bounds.   Specifically, by setting $h (x,t(x))$ in  Corollary~\ref{cor_7}  as 
\begin{align}\label{seq_5_eq_37}
 h (x,t(x)) =Q\big(x, t, M (t) \big),
\end{align}
where $M (t)$ is the MGF, and  $Q\big(x, t, M (t) \big)$ is any continuous function that satisfies $Q\big(x, t, M (t) \big)>0$ and   $\lim\limits_{x\to l}Q\big(x, t, M (t) \big)=0$,  and then inserting $h (x,t(x))$, given by \eqref{seq_5_eq_37}, into the optimization problems in \eqref{seq_5_nn1} and \eqref{seq_5_nn2}, we obtain corresponding upper and lower bounds that  depend on the MGF.

We have yet to provide functions $g (x)$ that are good candidates for the corresponding lower bound in Theorem~\ref{thm_2}. Such $g (x)$ will be arrived  at by an iterative method, which is the subject of the following section.

\section{The Iterative Method For The Left Tail}\label{sec-6}

In this section, we provide an iterative method for obtaining ever tighter upper and lower bounds on $F (x)$, under certain conditions. Before we provide the iterative method, we introduce several lemmas that will be useful  for arriving at the iterative method. Moreover, in this section, when we say that some function is a valid upper or lower bound on $F (x)$ as per Theorem~\ref{thm_2}, we mean that this bound is obtained using Theorem~\ref{thm_2} and thereby satisfies all of the conditions laid out in Theorem~\ref{thm_2}.

We start with the following lemma.
 \begin{lemma}\label{lema_6}
If $g (x)$, with properties defined in Theorem~\ref{thm_2}, also  satisfies the  limit  
 \begin{align}\label{seq_5_eq_13-n1}
 \lim_{x\to l} g (x)=0,
 \end{align}
 then $P (x)$, given by \eqref{seq_5_eq_7a},  satisfies \eqref{seq_5_eq_7b}, i.e.,  $P (x)$ satisfies the limit  
 \begin{align}\label{seq_5_eq_13-n2}
\lim_{x\to l} P (x)=  \lim_{x\to l}    f (x)\frac{g (x)}{g' (x)}=0.
\end{align}
 
 \end{lemma}
\begin{IEEEproof}
The proof is provided in Appendix~\ref{app_7}.
\end{IEEEproof}

We note that a function $g (x)$ does not need to satisfy   \eqref{seq_5_eq_13-n1} in order for \eqref{seq_5_eq_7b} (i.e.,  \eqref{seq_5_eq_13-n2}) to hold. In other words, there are functions $g (x)$ for which \eqref{seq_5_eq_13-n1} does not hold and yet \eqref{seq_5_eq_7b} (i.e.,  \eqref{seq_5_eq_13-n2}) holds. However, what Lemma~\ref{lema_6} shows  us is that if $g (x)$ is such that \eqref{seq_5_eq_13-n1} holds, then we have certainty that \eqref{seq_5_eq_7b} (i.e.,  \eqref{seq_5_eq_13-n2})   holds. We will find   Lemma~\ref{lema_6} useful later on.

We now start providing the basic building elements of the iterative method. 

Let us define $P_0(x)$ as
\begin{align}\label{seq_6_eq_1}
P_0(x)=   f (x)\frac{g (x)}{g' (x)}.
\end{align}
Note that $P_0(x)$ in \eqref{seq_6_eq_1} is identical to $P (x)$ given by \eqref{seq_5_eq_7a}.
Let us assume that $g (x)$ satisfies the conditions defined in Theorem~\ref{thm_2} in order for \eqref{seq_6_eq_1} to be an upper on $F (x)$, $\forall x<x_0$, or for  \eqref{seq_6_eq_1} to be a lower on $F (x)$, $\forall x<x_0$.
 Next, let us define the function $P_i(x)$, for $i=0,1,2,...$ as
 \begin{align}\label{seq_6_eq_2}
P_{i+1}(x)=  f (x)\frac{P_i(x)}{P'_i(x)}.
\end{align}
Note that $P_{i+1}(x)$ is obtained in an itterative manner starting from the seed $P_{0}(x)$ given by \eqref{seq_6_eq_1}.

For the limit of $P_{i+1}(x)$ as $x\to l$, we have the following lemma.

 \begin{lemma}\label{lema_6a}
If $P_0(x)$ is a  valid upper bound or a valid lower bound  on $F (x)$, $\forall x<x_0$, as per Theorem~\ref{thm_2}, then the following limit holds for any $i=0,1,2,...,$
 \begin{align}\label{seq_5_eq_13-n2-1}
\lim\limits_{x\to l} P_{i+1}(x)=0.
\end{align}
 \end{lemma}
 
 \begin{IEEEproof}
The proof is provided in Appendix~\ref{app_8}.
\end{IEEEproof}

Next, we have the following useful lemma for $P_{i+1}(x)$.

\begin{lemma}\label{lema_7}
If $P_{i}(x)$ is an upper on $F (x)$, $\forall x<x_0$, as per Theorem~\ref{thm_2}, which also satisfies $P'_{i}(x)>0$,  $\forall x<\hat x_0$, then the following holds
 \begin{align}\label{seq_6_eq_3}
P_{i+1}(x)\leq P_{i}(x), \; \forall x<x_i \textrm{ and any } i=0,1,2,... ,
\end{align}
where $x_i=\min(x_0,\hat x_0)$.

Otherwise, if $P_{i}(x)$ is a lower bound on $F (x)$, $\forall x<x_0$, as per Theorem~\ref{thm_2}, which also satisfies $P'_{i}(x)>0$,  $\forall x<\hat x_0$,    then the following holds
 \begin{align}\label{seq_6_eq_4}
P_{i+1}(x)\geq P_{i}(x),\; \forall x<x_i \textrm{ and any } i=0,1,2,...,
\end{align}
where $x_i=\min(x_0,\hat x_0)$.
\end{lemma}

 \begin{IEEEproof}
The proof is provided in Appendix~\ref{app_9}.
\end{IEEEproof}

Lemma~\ref{lema_7} is useful since it tells us that the next iteration  $P_{i+1}(x)$ of an upper bound $P_i(x)$ is always smaller than the preceding iteration  $P_i(x)$. Thereby, if the next iteration,  $P_{i+1}(x)$, itself is also an upper bound, then  $P_{i+1}(x)$ will be a tighter upper bound than its preceding iteration,  $P_i(x)$. Similarly, Lemma~\ref{lema_7}  tells us that the next iteration,  $P_{i+1}(x)$, of a lower bound $P_i(x)$ is always larger than the preceding iteration, $P_i(x)$. Thereby, if the new iteration, $P_{i+1}(x)$, itself is also a lower bound, then $P_{i+1}(x)$ will be a tighter lower bound than its preceding iteration, $P_{i}(x)$. Again, we can imagine that problems may arise if the next iteration of an upper bound (lower bound) becomes a lower bound (upper bound). But when that happens, we have methods to check if that will result in a tighter bound than the one from the previous iteration, as explained in the following.

We now introduce  an auxiliary lower bound which we can use to measure the tightness of a lower bound obtained iteratively from a preceding upper bound. Specifically,
 note that from a given upper bound, $P_i(x)$, we can always create an auxiliary  lower bound, denoted by $P_{L,i}(x)$, which is obtained by reflecting $P_i(x)$ with respect to $F (x)$, as
 \begin{align}\label{seq_6_eq_4a}
P_{L,i}(x)=  F (x)-\big(P_i(i)-F (x)\big) =  2 F (x) - P_i(x), \; \forall x<x_i,
\end{align}
when $P_i(x)$ is a valid upper  bound on $F (x)$  as per Theorem~\ref{thm_2}.

Similarly, we introduce  an auxiliary upper bound which we can use to measure the tightness of an upper bound obtained iteratively from a preceding lower bound. Again, note that from a given lower bound, $P_i(x)$, we can always create an auxiliary  upper bound, denoted by $P_{U,i}(x)$, which is obtained by reflecting $P_i(x)$ with respect to $F (x)$, as
 \begin{align}\label{seq_6_eq_4aa}
P_{U,i}(x)=  F (x)+\big(F (x)-P_i(i)\big) =  2F (x) - P_i(x), \; \forall x<x_i,
\end{align}
when $P_i(x)$ is a valid lower  bound on $F (x)$  as per Theorem~\ref{thm_2}.

Using $P_{L,i}(x)$ and $P_{U,i}(x)$, we can state the following lemma.

\begin{lemma}\label{lema_8}
Let $P_{i}(x)$ be an upper on $F (x)$, $\forall x<x_i$, and let
 \begin{align}\label{seq_6_eq_5}
P_{i+1}(x)=   f (x)\frac{P_i(x)}{P'_i(x)}, \; \forall x<x_{i+1},
\end{align}
be a lower bound on $F (x)$, $\forall x<x_{i+1}$. In that case, if the condition
\begin{align}\label{seq_6_eq_7}
P'_{i+1}(x)+P'_{i}(x)-2 f (x)\geq 0, \; \forall x<\hat x_{i+1},
\end{align}
holds, then 
$P_{i+1}(x)$ is a tighter bound on $F (x)$ than $P_{i}(x)$, in the sense that
 \begin{align}\label{seq_6_eq_8}
P_{L,i}(x)\leq  P_{i+1}(x), \; \forall x<\hat x_{i+1},
\end{align}
holds, where $P_{L,i}(x)$ is given by \eqref{seq_6_eq_4a}.

On the other hand, let $P_{i}(x)$ be a lower bound on $F (x)$, $\forall x<x_i$, and let
 \begin{align}\label{seq_6_eq_9}
P_{i+1}(x)=   f (x)\frac{P_i(x)}{P'_i(x)}, \; \forall x<x_{i+1},
\end{align}
be an upper bound on $F (x)$, $\forall x<x_{i+1}$. In that case, if the condition
 \begin{align}\label{seq_6_eq_11}
P'_{i+1}(x)+P'_{i}(x)-2 f (x)\leq 0, \; \forall x<\hat x_{i+1},
\end{align}
holds, then $P_{i+1}(x)$ is a tighter bound on $F (x)$ than $P_{i}(x)$, in the sense that
 \begin{align}\label{seq_6_eq_10}
P_{U,i}(x)\geq  P_{i+1}(x), \; \forall x<\hat x_{i+1},
\end{align}
 where $P_{U,i}(x)$ is given by \eqref{seq_6_eq_4aa}.

\end{lemma}

 \begin{IEEEproof}
The proof is provided in Appendix~\ref{app_10}.
\end{IEEEproof}

We now have all of the necessary elements to provide an iterative algorithm that can lead to ever tighter upper and lower bounds obtained in an iterative manner; bounds which not necessary all  hold up to the same $x=x_0$.  Specifically, if $P_{i}(x)$ is an upper/lower bound on $F (x)$, $\forall x<x_i$, then it may happen that $P_{i+1}(x)$ is a lower/upper bound on $F (x)$, $\forall x<x_{i+1}$, where $x_{i+1}<x_i$.

\begin{algorithm}
\caption{The Iterative Method For The Left Tail}\label{alg:2}
 
\SetKwData{Left}{left}
\SetKwInOut{Set}{Set}
	\SetKwInOut{Ensure}{Ensure}
	\SetKwInOut{Input}{Input}
	\SetKwInOut{Output}{Output}
	\SetKwInOut{Return}{Return}
	\SetKwData{And}{and}
	\Input{$f (x)$, $g (x)$, $x_0$}
	\Output{$P_L(x)$, $P_U(x)$}
	\BlankLine
	\Set{$P_0(x)=   f (x)\dfrac{g (x)}{g' (x)}$}
	\BlankLine
	\Set{$C=0$, $i=0$, $P_L(x)=NaN$, $P_U(x)=NaN$}
	\While{$C=0$  }{
	$P_{i+1}(x)=   f (x)\dfrac{P_{i}(x)}{P'_{i}(x)}$\\
	\uIf{$P_{i}(x)> 0, \; \forall x<x_0$, \textbf{ and } $P'_{i}(x)> 0, \; \forall x<x_0$}
	{
	\uIf{$P'_{i}(x)-f (x)\geq 0, \; \forall x<x_0$,}
	{
	\uIf{$P'_{i+1}(x)-f (x)\geq 0, \; \forall x<x_0$,}
	{	\BlankLine
	$P_U(x)=P_{i+1}(x)$
		\BlankLine}
	\uElseIf{$P'_{i+1}(x)-f (x)\leq 0, \; \forall x<x_0$, \textbf{ and } $P'_{i+1}(x)+P'_{i}(x)-2f (x)\geq 0, \; \forall x<x_0$,}
	{	\BlankLine
	$P_L(x)=P_{i+1}(x)$
		\BlankLine}
	\Else{$C=1$}
	}
	\Else{
	\uIf{$P'_{i+1}(x)-f (x)\geq 0, \; \forall x<x_0$, \textbf{ and } $P'_{i+1}(x)+P'_{i}(x)-2f (x)\leq 0, \; \forall x<x_0$,}
	{$P_U(x)=P_{i+1}(x)$}
	\uElseIf{$P'_{i+1}(x)-f (x)\leq 0, \; \forall x<x_0$,}
	{$P_L(x)=P_{i+1}(x)$}
	\Else{$C=1$}
	}
	}
	\Else{$C=1$}
		$i=i+1$
	}
	\Return{$P_L(x)$, $P_R(x)$}
\end{algorithm}

The iterative algorithm is given in Algorithm~\ref{alg:2}, and works as follows.  The algorithm takes as inputs the PDF, $f (x)$, the function $g (x)$, and a desired point $x_0$ up to which we want these iteratively obtained bounds to hold. Note, the function $g (x)$ must be such that $P_0(x)$ is a valid upper or lower bound $\forall x<x_0$, as per Theorem~\ref{thm_2}. The \textit{while loop}  performs the following computations in an iterative manner, unless in the process of iteration  $C$ changes value from $C=0$ to $C=1$.
The   \textit{outer if  condition}, checks if the function $g (x)$ in iteration $i$, which in this case is $g (x)=P_i(x)$, is a valid function according to Theorem~\ref{thm_2}. If true, then the algorithm continues to the \textit{middle if  condition}. If false, $C$ is set to $C=1$ and the  \textit{while loop} stops.
The   \textit{middle if  condition} in the algorithm checks  if the previous iteration $P_i(x)$ is a valid upper bound. If true, then the \textit{inner if  condition}   checks if the next iteration $P_{i+1}(x)$ is also a valid  upper bound. If true, then $P_{i+1}(x)$ is a tighter upper bound than $P_i(x)$, according to Lemma~\ref{lema_7},  and therefore   $P_{i+1}(x)$ is stored into $P_U(x)$. Otherwise, if the \textit{inner if  condition} is false,  it is checked whether the next iteration $P_{i+1}(x)$ is  a valid  lower bound and if $P_{i+1}(x)$ is a tighter bound  than $P_{i}(x)$, as per Lemma~\ref{lema_8}. If true, then   $P_{i+1}(x)$ is stored into $P_L(x)$. If false, $C$ is set to $C=1$ and the \textit{while loop} stops.
 On the other hand, if the \textit{middle if  condition} is false, then the previous iteration $P_i(x)$ must be a valid lower bound. Therefore,  the \textit{inner if  condition} checks if the next iteration $P_{i+1}(x)$ is   a valid  upper bound and if this bound is tighter than $P_{i}(x)$, as per Lemma~\ref{lema_8}. If true, then $P_{i+1}(x)$ is stored into $P_U(x)$. If false,  then the next iteration $P_{i+1}(x)$ is checked whether it is valid  lower bound. If true, then $P_{i+1}(x)$ is a tighter lower bound than $P_i(x)$, according to Lemma~\ref{lema_7},  and therefore   $P_{i+1}(x)$  is stored   into $P_L(x)$. If false,    $C$ is set to $C=1$ and the \textit{while loop} stops. If $C$ has not changed value during one cycle in the \textit{while loop}, then $C$  keeps the value  $C=0$ and therefore the \textit{while loop} performs another iteration. 
 Finally, when the \textit{while loop} stops, the algorithm returns $P_L(x)$ and $P_U(x)$.

%%%%%%%%%%%%%%%%%%%%%%%%%%%%%%%%%%%%%%%%%%%%%%%%%%%%%%%%%%%%%%%%%%%%%%%%%%%%%%%%%%%%%%%%%%%%%%%%%%%%%%%%%%%%%%%%%%%%%%%%%%%%%%%%%%%%%%%%%%%%%%%%

\section{Rate of Convergence}\label{sec-7}

Since, in general, we can obtain upper and  lower bounds on $1-F (x)$ and on $F (x)$, we can measure how fast an upper bound and a lower bound  converge to each other using the rate of convergence function, given by
\begin{align}\label{sec_8_eq_1}
R (x)=\frac{P_U(x)}{P_L(x)}-1 =  \frac{g_U(x) g'_L(x)}{g_L(x) g'_U(x)}-1, 
\end{align}
where $P_U(x)$ and $P_L(x)$ are upper and lower bounds on $1-F (x)$, constructed using $g_U(x)$ and $g_L(x)$, respectively, as per Theorem~\ref{thm_1}, or  $P_U(x)$ and $P_L(x)$ are upper and lower bounds on $F (x)$, constructed using $g_U(x)$ and $g_L(x)$, respectively, as per Theorem~\ref{thm_2}.

\begin{remark}
We note that expression \eqref{sec_8_eq_1} might also be helpful towards the search for the optimal functions $g_U(x)$ and $g_L(x)$ that result in the tightest upper and lower bounds. Specifically,  the optimal functions $g_U(x)$ and $g_L(x)$  are the ones that minimize \eqref{sec_8_eq_1} under the constraint that $g_U(x)$ and $g_L(x)$ satisfy the conditions for $g (x)$ laid out in Theorem~\ref{thm_1} or Theorem~\ref{thm_2}, and $g_U(x)$ and $g_Lx)$  are given in the form of closed-form expressions.
\end{remark}

If we use the iterative method to obtain a lower bound, $P_{i+1}$, from an upper bound, $P_{i}$, on $1-F (x)$, as per Algorithm~\ref{alg:1}, then the rate of convergence would be
\begin{align}\label{sec_8_eq_2}
R (x)= \frac{P_{i}(x)}{P_{i+1}(x)}-1 = \frac{P_{i}(x)}{-f (x)\dfrac{P_{i}(x)}{P'_{i}(x)}}-1  = \frac{-P'_{i}(x)}{f (x)}-1.
\end{align}

Similarly,  if we use the iterative method to obtain an upper bound, $P_{i+1}$, from a lower bound, $P_{i}$, on $1-F (x)$, as per Algorithm~\ref{alg:1}, then the rate of convergence would be
\begin{align}\label{sec_8_eq_3}
R (x)= \frac{P_{i+1}(x)}{P_{i}(x)}-1 = \frac{-f (x)\dfrac{P_{i}(x)}{P'_{i}(x)}}{P_{i}(x)}-1  = \frac{f (x)}{-P'_{i}(x)}-1.
\end{align}

On the other hand, if we use the iterative method to obtain a lower bound, $P_{i+1}$, from a upper bound, $P_{i}$, on $F (x)$, as per Algorithm~\ref{alg:2}, then the rate of convergence would be
\begin{align}\label{sec_8_eq_4}
R (x)= \frac{P_{i}(x)}{P_{i+1}(x)}-1 = \frac{P_{i}(x)}{f (x)\dfrac{P_{i}(x)}{P'_{i}(x)}}-1  = \frac{P'_{i}(x)}{f (x)}-1.
\end{align}

Similarly,  if we use the iterative method to obtain an upper bound, $P_{i+1}$, from a lower bound, $P_{i}$, on $F (x)$, as per Algorithm~\ref{alg:2}, then the rate of convergence would be
\begin{align}\label{sec_8_eq_5}
R (x)= \frac{P_{i+1}(x)}{P_{i}(x)}-1 = \frac{f (x)\dfrac{P_{i}(x)}{P'_{i}(x)}}{P_{i}(x)}-1  = \frac{f (x)}{P'_{i}(x)}-1.
\end{align}

Note that the rate of convergence can always be obtained in a closed-form expression, given that $g_U(x)$ and $g_L(x)$, i.e., $P'_{i}(x)$, are given  in a closed-form expression.  The rate of convergence is important since it provides information on how close the bounds are to $1-F (x)$ or to $F (x)$, without having any information about $1-F (x)$ or $F (x)$.

\section{Bounding and Approximating The Converse Bound Of The AWGN Channel Capacity}\label{sec-awgn}

One application of the proposed tail-bounding methods is bounding and approximating the converse  bound of the AWGN channel capacity in the finite-blocklength regime, as shown in this section.

\subsection{Problem Formulation}
Assume an AWGN channel with SNR $\Omega$. Assume transmission of codewords of length $n$ on this channel. Assume a prescribed average error rate of the codewords, denoted by $\epsilon$. Then, as proved in \cite{5452208},  the capacity of this channel $C(n,\epsilon)$ is upper bounded by  $R(n,\epsilon)$,  given by
\begin{align}\label{eq:cap_1}
C(n,\epsilon) \leq  R(n,\epsilon)= -\frac{1}{n} \log_2\left(F_{\rm FA}\left(\frac{n\lambda}{1+\Omega}\right) \right),
\end{align}
where $\lambda$ is found from 
\begin{align}\label{eq:cap_2}
1-F_{\rm MD}(n\lambda)=\epsilon,  
\end{align}
where $F_{\rm FA}(x)$ is the CDF of the non-central chi-squared RV with degrees of freedom $n$ and non-centrality parameter $n\frac{1+\Omega}{\Omega}$ and $F_{\rm MD}(x)$ is the CDF of the non-central chi-squared RV with degrees of freedom $n$ and non-centrality parameter $\frac{n}{\Omega}$, see \cite{7303958}. The CDF of the  non-central  chi-squared RV with $n$ degrees-of-freedom and non-centrality parameter $s$ is given by
\begin{align}\label{eq:marqum_Q}
F(x)=1-Q_{\frac{n}{2}}\left(\sqrt{s },\sqrt{x}\right),
\end{align}
where $Q_{M}(a,b)$ is the  Marcum-Q function.
Of course, evaluating the capacity bounds in \eqref{eq:cap_1} is difficult due to the numerical instability  of the Marcum-Q function in \eqref{eq:marqum_Q}, and as a result efforts have been made in \cite{7303958} and \cite{7589108} to obtain a more numerically pleasant expressions for numerical evaluation of \eqref{eq:cap_1}.  

\subsection{Lower and Upper Bounds On The Converse Bound}
Using the proposed tail-bounding method, we can derive closed-form and tight upper and lower bounds on the left tail of $F_{\rm FA}\left(\frac{n\lambda}{1+\Omega}\right) $ in \eqref{eq:cap_1}, and thereby obtain respective closed-form and tight lower and upper bounds on the converse $R(n,\epsilon)$. Next, we can derive closed-form and tight upper and/or lower bounds on the right-tail of $1-F_{\rm MD}(n\lambda) $ in \eqref{eq:cap_2} and thereby find $\lambda$ using a numerically stable method even for very large $n$ and $\Omega$. 

For the right tail, $1-F_{\rm MD}\left( n\lambda \right) $, we   use the proposed tail-bounding method  to obtain  an upper and a lower bound using  $g (x)= f (x)$ (and more tighter bound can be obtained using the iterative method, if desired). Thereby, we obtain the following   expressions, which holds for $n>n_0$,
\begin{align}
&1-F_{\rm MD}\left( n\lambda \right)\leq P_{0,{\rm MD}} \left( n\lambda \right)
\nonumber\\
&=  
\frac{-e^{-\frac{n(1+\lambda\Omega)}{2\Omega}}\;
n\sqrt{\lambda/\Omega}\;(\lambda\Omega)^{n/4}\;
I_{\frac{n}{2}-1}\!\bigl(n\sqrt{\lambda/\Omega}\bigr)^2}{\;
(-2+n-n\lambda) I_{\frac{n}{2}-1}\!\bigl(n\sqrt{\lambda/\Omega}\bigr)
+ n\sqrt{\lambda/\Omega} I_{\frac{n}{2}}\!\bigl(n\sqrt{\lambda/\Omega}\bigr)}
\label{eq:F_MA_bound1}\\
  &1-F_{\rm MD}\left( n\lambda \right)\geq P_{1,{\rm MD}} \left( n\lambda \right)
\nonumber\\
&=  
f\left(n\lambda\right)  \frac{P_{0,{\rm MD}}\left(n\lambda\right) }{P_{0,{\rm MD}}'\left(n\lambda\right) } ,
\label{eq:F_MA_bound2}  
\end{align} 
where $n_0$ can be found as the integer $n$ for which the bound goes beyond one.
Now inserting   $P_{0,{\rm MD}} \left( n\lambda \right)$  instead of $1-F_{\rm MD}(n\lambda)$ in \eqref{eq:cap_2}, we obtain
\begin{align}\label{eq:cap_2-1}
P_{0,{\rm MD}} \left( n\lambda \right) =\epsilon,  
\end{align}
and
\begin{align}\label{eq:cap_2-2}
P_{1,{\rm MD}} \left( n\lambda \right) =\epsilon,  
\end{align}
which are much easier to solve numerically for $\lambda$ than \eqref{eq:cap_2}.

We  now use the proposed iterative method  to first obtain   lower and upper bounds of $F_{\rm FA}\left(\frac{n\lambda}{1+\Omega}\right) $, utilizing $g (x)=x f (x)$ for the non-central chi-squared RV. Then, we plug these bounds into \eqref{eq:cap_1} instead of $F_{\rm FA}\left(\frac{n\lambda}{1+\Omega}\right) $, and obtain the following lower and upper bounds on the converse bounds  
\begin{align}
R(n,\epsilon) &\geq -\frac{1}{n} \log_2\left( P_{0,{\rm FA}}\left(\frac{n\lambda}{1+\Omega}\right)  \right)
 \label{eq:F_FA_bound1}\\
R(n,\epsilon)  & \leq -\frac{1}{n} \log_2\left(  P_{1,{\rm FA}}\left(\frac{n\lambda}{1+\Omega}\right)\right) \label{eq:F_FA_bound2}, 
\end{align} 
where   $\lambda$ in \eqref{eq:F_FA_bound1} and in \eqref{eq:F_FA_bound2} are obtained from \eqref{eq:cap_2-1} and \eqref{eq:cap_2-2}, respectively, and where
\begin{align}
&P_{0,{\rm FA}}\left(\frac{n\lambda}{1+\Omega}\right) 
\nonumber\\
&= \frac{
  \E^{-\frac{n(2+\lambda+1/\Omega+\Omega)}{2(1+\Omega)}}
  \sqrt{\lambda}\big(\tfrac{\lambda\Omega}{(1+\Omega)^2}\big)^{\!n/4}(1+\Omega)
  \I_{\frac{n-2}{2}}\!\big(n\sqrt{\tfrac{\lambda}{\Omega}}\big)^{\!2}
}{
  \sqrt{\Omega}(1-\lambda+\Omega)\I_{\frac{n-2}{2}}\!\big(n\sqrt{\tfrac{\lambda}{\Omega}}\big)\;+\;
  \sqrt{\lambda}(1+\Omega)\I_{n/2}\!\big(n\sqrt{\tfrac{\lambda}{\Omega}}\big)
}, \label{eq:F_FA_bound1a}\\
&  P_{1,{\rm FA}}\left(\frac{n\lambda}{1+\Omega}\right) =f\left(\frac{n\lambda}{1+\Omega}\right)  \frac{P_{0,{\rm FA}}\left(\frac{n\lambda}{1+\Omega}\right) }{P_{0,{\rm FA}}'\left(\frac{n\lambda}{1+\Omega}\right) } \label{eq:F_FA_bound2a}.
\end{align}

\subsection{Closed-Form Asymptotic Approximation Of  The Converse Bound}

We now turn on deriving a closed-form asymptotic approximation of  the converse bound in \eqref{eq:cap_1}. To this end, we first find an asymptotic  solution for  $\lambda$ from \eqref{eq:cap_2-1}, which  holds for large $n$, which is obtained as
\begin{align}
\lambda \;&\approx\; 
1+\frac{1}{\Omega}\;+\;\frac{1}{\sqrt{n}}\,
\sqrt{\frac{2(\Omega+2)}{\Omega}\;
W\!\left(\frac{1}{2\pi\epsilon^2}\right)}.
\label{eq:final11}
\end{align}
 The full derivation of \eqref{eq:final11} is shown in Appendix~\ref{sec-app-lam}. 

Next, we simplify \eqref{eq:F_FA_bound1}   for large $n$ and thereby we obtain
 \begin{align}\label{eq:My_rate1}
  R(n, \epsilon)  &\approx \frac{1}{2}\log_2\!\big(1+\Omega \big) -\frac{1}{2}\log_2\!\left(\frac{2\lambda}{ 1+\sqrt{1+\frac{4\lambda}{\Omega}}}\right) \nonumber\\
 &+\frac{1}{2\ln(2)}\left(1+\frac{1}{\Omega}+\frac{\lambda}{1+\Omega}- \sqrt{1+\frac{4\lambda}{\Omega}}\right),
 \end{align}
 where $\lambda$ is given by \eqref{eq:final11}.
 The   derivation of \eqref{eq:My_rate1} is shown in Appendix~\ref{sec-app-R}. 

We can observe from \eqref{eq:final11} that  $\lim\limits_{n\to\infty}\lambda=1+\frac{1}{\Omega}$, if $\epsilon>0$. Plugging in  $\lambda=1+\frac{1}{\Omega}$ into \eqref{eq:My_rate1} we obtain that $R(n, \epsilon)= \frac{1}{2}\log_2 (1+\Omega )$. Thereby, for $\epsilon>0$, the  following holds 
 \begin{align}
 \lim\limits_{n\to \infty}R(n, \epsilon)=\frac{1}{2}\log_2\!\big(1+\Omega \big),
 \end{align}
 which is the channel capacity of the infinite blocklength. Hence, indeed, the derived expression is asymptotically correct. We will see   via numerical examples in the next section that \eqref{eq:My_rate1} is much tighter approximation to the converse bound than \eqref{eq:Polyanskiy_rate}.

\section{Numerical Examples}\label{sec-num}

In this section, we apply the proposed method for upper and lower bounding the right tail of the  Gaussian and the beta prime RVs, whose CDFs are given by
\begin{align}
F (x)&=\frac{1}{2} \left(\mathrm{erf}\left(\frac{x-\mu }{\sqrt{2} \sigma }\right)+1\right),\\
F (x)&=\frac{B_{\frac{x}{x+1}}(\alpha ,\beta )}{B(\alpha ,\beta )},
\end{align}
respectively, where $\mathrm{erf}(x)$, $B_x(a,b)$, and $B(a,b)$ are the Gaussian error function, the incomplete Beta function, and the Beta function, respectively. 

In addition, we apply the proposed method for upper and lower bounding the left tail of the   non-central  chi-squared RV, whose CDF is given by
\begin{align}
F (x)=1-Q_{\frac{k}{2}}\left(\sqrt{\lambda },\sqrt{x}\right),
\end{align}
where $Q_{M}(a,b)$ is the  Marcum-Q function.

Finally, as a third example, we use the proposed tail-bounding method to derive the approximate finite-blocklength data rate

For both the right and the left tails, we will use the iterative method to arrive at ever tighter upper and lower bounds on the tails, starting from a seed.

Since the bounds that we will illustrate are very tight, for a better visual representation, we  may choose to plot the functions
\begin{align}\label{sec-num-eq_1}
\left|\frac{1-F (x)}{P_{i}(x)} - 1\right|
\end{align}
and
\begin{align}\label{sec-num-eq_1a}
\left|\frac{F (x)}{P_{i}(x)} - 1\right|
\end{align}
for the right and the left tail bounds, respectively, for different $i$'s, where $F (x)$ is evaluated numerically. The functions in  \eqref{sec-num-eq_1} and \eqref{sec-num-eq_1a}  show  how fast the bound $P_{i}(x)$ converges to the right tail or to the left tail, respectively, independent of whether $P_{i}(x)$ is an upper or lower bound on the corresponding tail.
 However, continuing with the assumption that we do not have any access to $F (x)$, we choose instead to plot the rate of convergence, defined in Sec.~\ref{sec-7}, which provides information on how fast the upper and  lower bounds converge to each other. The rate of convergence can be written as
\begin{align}\label{sec-num-eq_2}
R_{i}(x)= \left|\frac{P_{i+1}(x)}{P_{i}(x)}-1\right|,
\end{align}
in both cases when $P_{i+1}(x)$ and $P_{i}(x)$ are upper and lower bounds, respectively, and when $P_{i+1}(x)$ and $P_{i}(x)$ are lower and upper bounds, respectively. Note that
\begin{align}
\left|\frac{1-F (x)}{P_{i}(x)} - 1\right| &\leq \left|\frac{P_{i+1}(x)}{P_{i}(x)}-1\right|\label{sec-num-eq_3a}
\end{align}
and
\begin{align}
\left|\frac{F (x)}{P_{i}(x)} - 1\right| &\leq \left|\frac{P_{i+1}(x)}{P_{i}(x)}-1\right|\label{sec-num-eq_3b}
\end{align}
always hold for the right and the left tail bounds, respectively. Hence,  from \eqref{sec-num-eq_3a} and \eqref{sec-num-eq_3b}, we see that the rate of convergence, $R_{i}(x)$, given by \eqref{sec-num-eq_2},  also provides an upper bound on the convergence of the bound $P_{i}(x)$ towards $1-F (x)$ or $F (x)$, independent of whether  $P_{i}(x)$ is an upper or lower bound on $1-F (x)$ or on $F (x)$.

\subsection{The Right Tail}

For the right tail, we  will use the iterative method to obtain  ever tighter closed-form lower and upper bounds. To this end, for constructing the seed, $P_0(x)$, we will use $g (x)=f (x)$ for the Gaussian RV and we will use $g (x)=(x-l)f (x)$ for the beta prime RV. We note that choosing $g (x)=f (x)$ is appropriate for RVs whose tail decays exponentially, whereas $g (x)=(x-l)f (x)$ is appropriate for RVs whose tail decays sub-exponentially.

\begin{figure}
\includegraphics[width=1\textwidth]{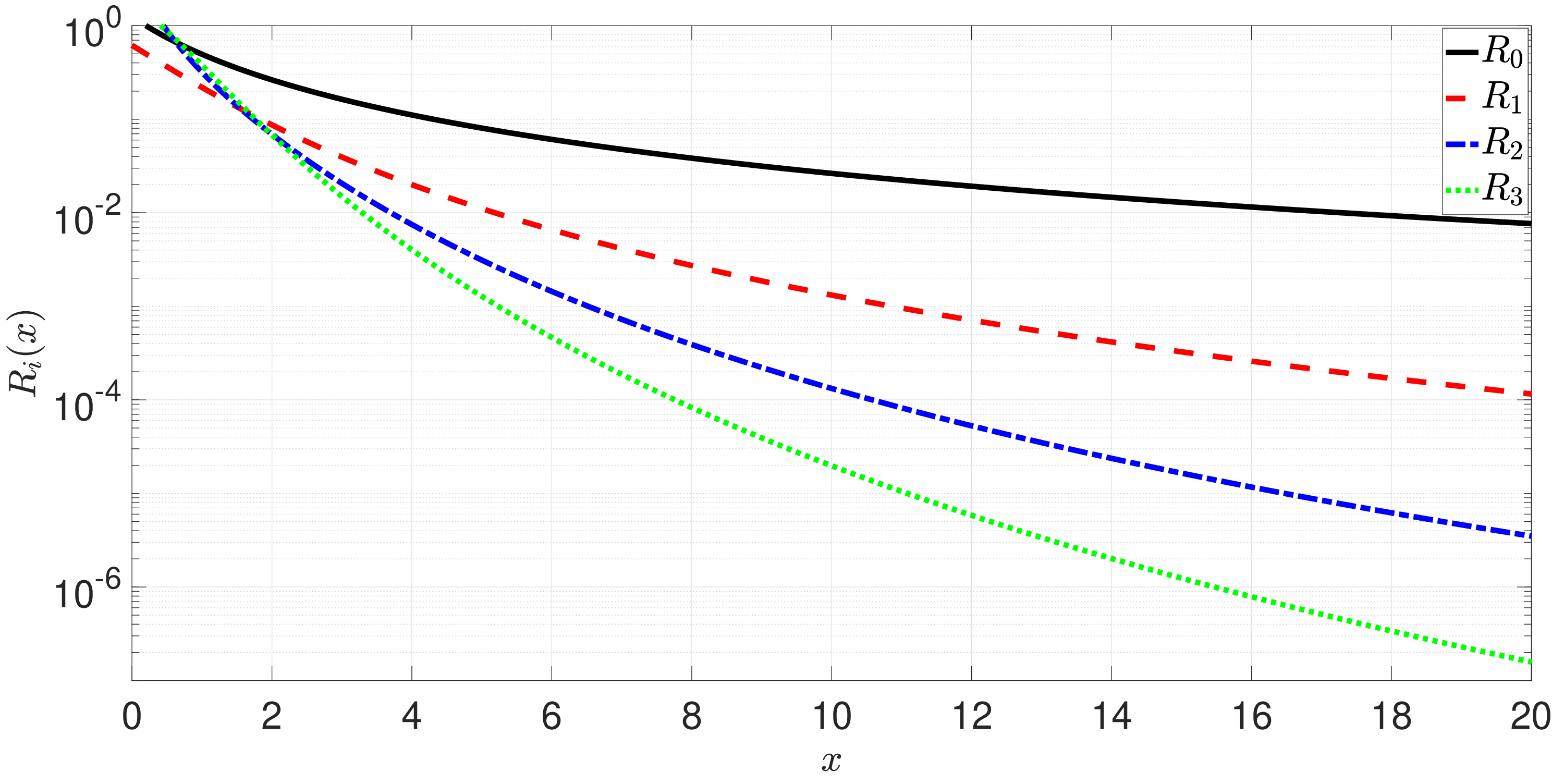}
\caption{Rate of convergence for the upper and lower bounds on  the right tail   of   the Gaussian RV with $\sigma=1.9$, $\mu=-1.7$.}
\label{fig_1}
\end{figure}

For the Gaussian right tail bounds, using the function $g (x)=f (x)$ to construct  the seed, we obtain the following closed-form  expressions for the bounds 
\begin{align}
1-F (x)&\leq P_0(x)=\frac{\sigma  e^{-\frac{(x-\mu )^2}{2 \sigma ^2}}}{\sqrt{2 \pi } (x-\mu )}, \;\forall x>\mu\\
1-F (x)&\geq P_1(x) =\frac{\sigma  (x-\mu ) e^{-\frac{(x-\mu )^2}{2 \sigma ^2}}}{\sqrt{2 \pi } \left(\sigma ^2+(x-\mu )^2\right)}, \;\forall x>\mu\\
1-F (x)&\leq P_2(x)
\nonumber\\
&\hspace{-4mm}=\frac{\sigma  (x-\mu ) e^{-\frac{(x-\mu )^2}{2 \sigma ^2}} \left(\sigma ^2+(x-\mu )^2\right)}{\sqrt{2 \pi } \left(-\sigma ^4+2 \sigma ^2 (x-\mu )^2+(x-\mu )^4\right)},\;\forall x>x_2\\
1-F (x)&\geq P_3(x)=e^{-\frac{(x-\mu )^2}{2 \sigma ^2}}
\nonumber\\
&\hspace{-14mm}\times  \frac{ \left(\sigma ^7 (\mu -x)+\sigma ^5 (x-\mu )^3+3 \sigma ^3 (x-\mu )^5+\sigma  (x-\mu )^7\right)}{\sqrt{2 \pi } \left(\sigma ^4+(x-\mu )^4\right)
   \left(\sigma ^4+4 \sigma ^2 (x-\mu )^2+(x-\mu )^4\right)},\nonumber\\
   & \hspace{-14mm} \;\forall x>x_3,\\
   1-F (x)&\leq P_4(x) ,\;\forall x>x_4,
\end{align} 
where the expressions for $x_2$, $x_3$, and   $x_4$ are omitted since they are too large to be fit in one row. Instead, they can be easily visualized from Fig.~\ref{fig_1}.
We also have a closed-form expression for $P_4(x)$, which we use  to plot $R_3(x)$ in Fig.~\ref{fig_1}, but we omit to show it analytically since it is too large to be fit in one row.

In Fig.~\ref{fig_1}, we show the rate of convergence, given by \eqref{sec-num-eq_2}, for $i=0,\; 1,\; 2,$ and $3$,  of the upper and lower bounds on  the right tail of   the Gaussian RV with $\sigma=1.9$ and $\mu=-1.7$.
Fig.~\ref{fig_1} shows that the upper and lower bounds on the Gaussian tail converge to each other very fast. Specifically, $P_0(x)$ and $P_1(x)$ converge to each other with rate proportional to $1/x^2$,  $P_1(x)$ and $P_2(x)$ converge to each other with rate proportional to $1/x^4$, $P_2(x)$ and $P_3(x)$ converge to each other with rate proportional to $1/x^6$, and $P_3(x)$ and $P_4(x)$ converge to each other with rate proportional to $1/x^8$.

 \begin{figure}
\includegraphics[width=1\textwidth]{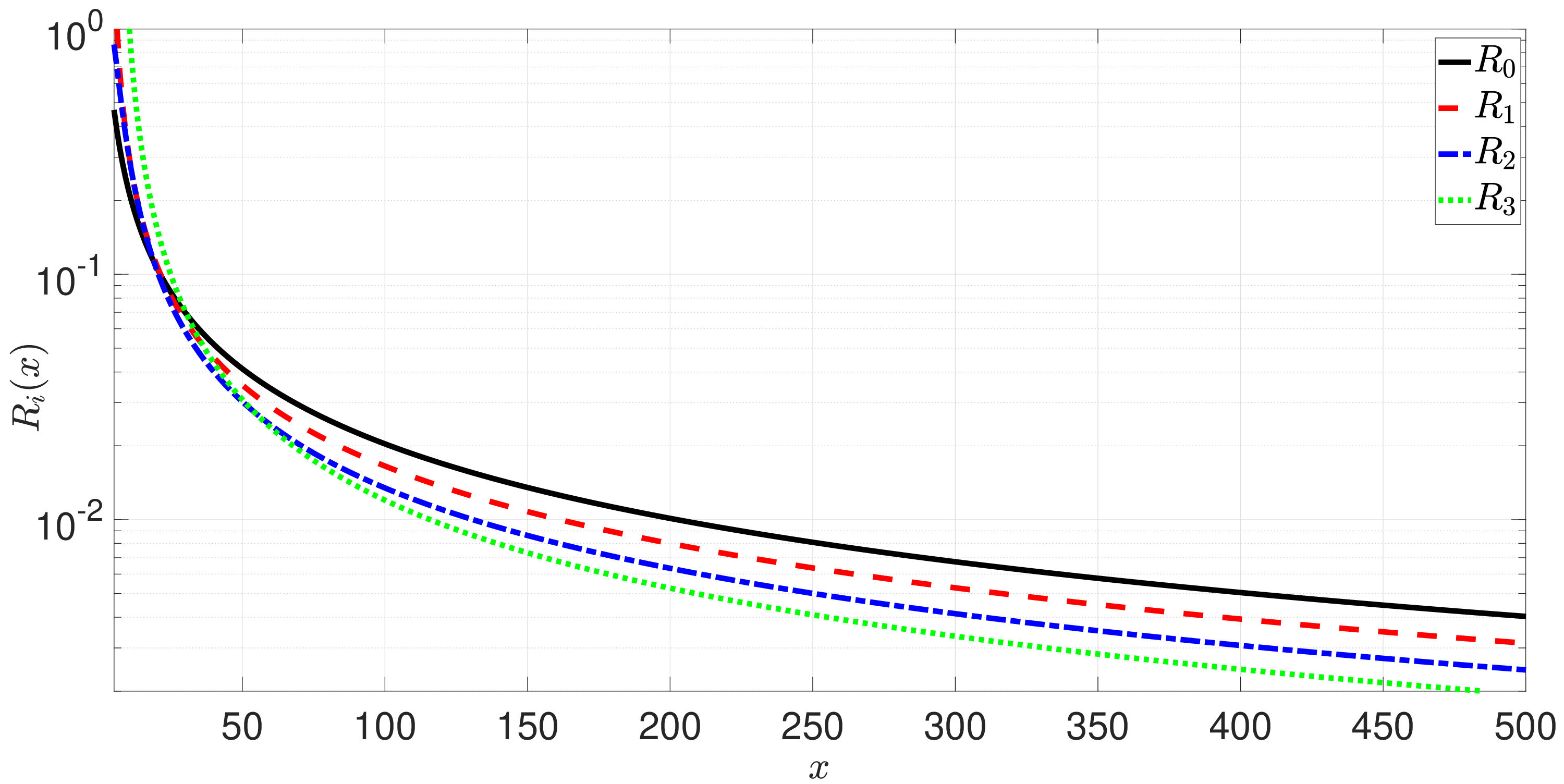}
\caption{Rate of convergence for the upper and lower bounds on   the right tail   of   the beta prime RV with $\alpha=2.1$, $\beta=1.3$.}
\label{fig_2}
\end{figure}

For the beta prime right tail bounds, using the function $g (x)=x f (x)$ (note $l=0$) to construct  the seed, we obtain the following   expressions for the bounds
\begin{align}
1-F (x)&\leq P_0(x) =\frac{x^{\alpha } (x+1)^{-\alpha -\beta +1}}{B(\alpha ,\beta ) (\beta  x-\alpha )}, \;\forall x>\frac{\alpha}{\beta}\\
1-F (x)&\geq P_1(x)
\nonumber\\
&\hspace{-12mm}
 =\frac{x^{\alpha } (x+1)^{-\alpha -\beta +1} (\beta  x-\alpha )}{B(\alpha ,\beta ) \left(\alpha ^2+\beta ^2 x^2+x (-2 \alpha  \beta +\alpha +\beta
   )\right)} , \;\forall x>\frac{\alpha}{\beta},\\
   1-F (x)&\leq P_2(x) ,\\
   1-F (x)&\geq P_3(x) ,\\
   1-F (x)&\leq P_4(x) . 
\end{align} 
We also have  expressions for $P_2(x)$, $P_3(x)$, and  $P_4(x)$, which we use  to plot $R_1(x)$, $R_2(x)$, and $R_3(x)$ in Fig.~\ref{fig_2}, but we omit to show them analytically since they are too large to be fit within one row.

In Fig.~\ref{fig_2}, we show the rate of convergence, given by \eqref{sec-num-eq_2}, for $i=0,\; 1,\; 2,$ and $3$, of the upper and lower bounds on  the right tail of the beta prime RV with $\alpha=2.1$ and  $\beta=1.3$.
Fig.~\ref{fig_2} shows that the upper and lower bounds on the right tail of the beta prime RV, specifically, the pairs $P_0(x)$ and $P_1(x)$,  $P_1(x)$ and $P_2(x)$,  $P_2(x)$ and $P_3(x)$, and $P_3(x)$ and $P_4(x)$ converge to each other with rate proportional to $1/x$.

\subsection{The Left  Tail}

For the left tail, we  also use the iterative method to obtain  ever tighter  lower and upper bounds. To this end, for constructing the seed, $P_0(x)$, we will use $g (x)=(x-l) f (x)$ for the non-central chi-squared RV (note $l=0$).

 \begin{figure}
\includegraphics[width=1\textwidth]{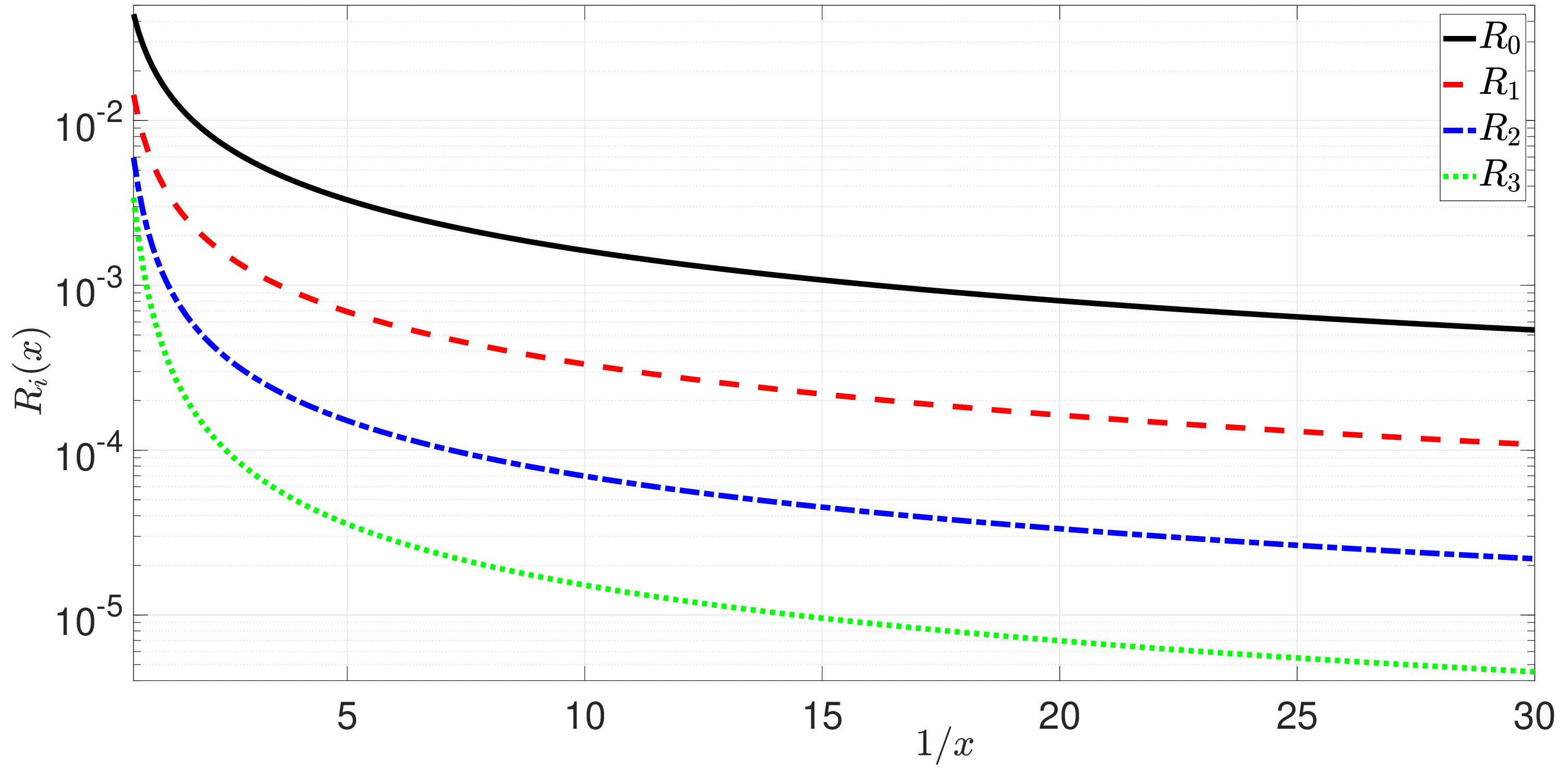}
\caption{Rate of convergence for the upper and lower bounds on the  left tail   for   the non-central chi-squared RV with $k=10$, $\lambda=2$.}
\label{fig_3}
\end{figure}

For the left tail bounds of the non-central chi-squared RV, using the function $g (x)=x f (x)$ to construct  the seed, we obtain the following   expressions 
\begin{align}
F (x)&\leq P_0(x)
\nonumber\\
&\hspace{-12mm} = \frac{e^{\frac{1}{2} (-\lambda -x)} \sqrt{\lambda  x} \left(\frac{x}{\lambda }\right)^{k/4} I_{\frac{k-2}{2}}\left(\sqrt{x \lambda }\right){}^2}{(k-x-2) I_{\frac{k-2}{2}}\left(\sqrt{x
   \lambda }\right)+\sqrt{\lambda  x} I_{\frac{k}{2}}\left(\sqrt{x \lambda }\right)},\; \forall x<x_0, \qquad \label{eq_p0_1}\\
   F (x)&\geq P_1(x) ,\\
F (x)&\leq P_2(x) ,\\
   F (x)&\geq P_3(x) ,\\
   F (x)&\leq P_4(x) .    
\end{align} 
where $x_0$ can be found as the point $x$ for which the bound goes beyond one. We also have   expressions for $P_1(x)$, $P_2(x)$, $P_3(x)$, and  $P_4(x)$, which we use  to plot $R_1(x)$, $R_2(x)$, and $R_3(x)$ in Fig.~\ref{fig_3}, but we omit to show them analytically since they are too large to be fit in one row.

In Fig.~\ref{fig_3}, we show the rate of convergence, given by \eqref{sec-num-eq_2}, for $i=0,\; 1,\; 2,$ and $3$, of the upper and lower bounds on  the left tail of the  non-central chi-squared RV with  $k=10$ and $\lambda=2$.
Fig.~\ref{fig_3}  shows that the upper and lower bounds on the left tail of the non-central chi-squared RV, specifically, the pairs $P_0(x)$ and $P_1(x)$,    $P_1(x)$ and $P_2(x)$,  $P_2(x)$ and $P_3(x)$, and $P_3(x)$ and $P_4(x)$ converge to each other with rate proportional to $x$, as $x\to 0$.

\subsection{Application: The Finite Blocklength Channel Capacity}

\begin{figure}
\includegraphics[width=1\textwidth]{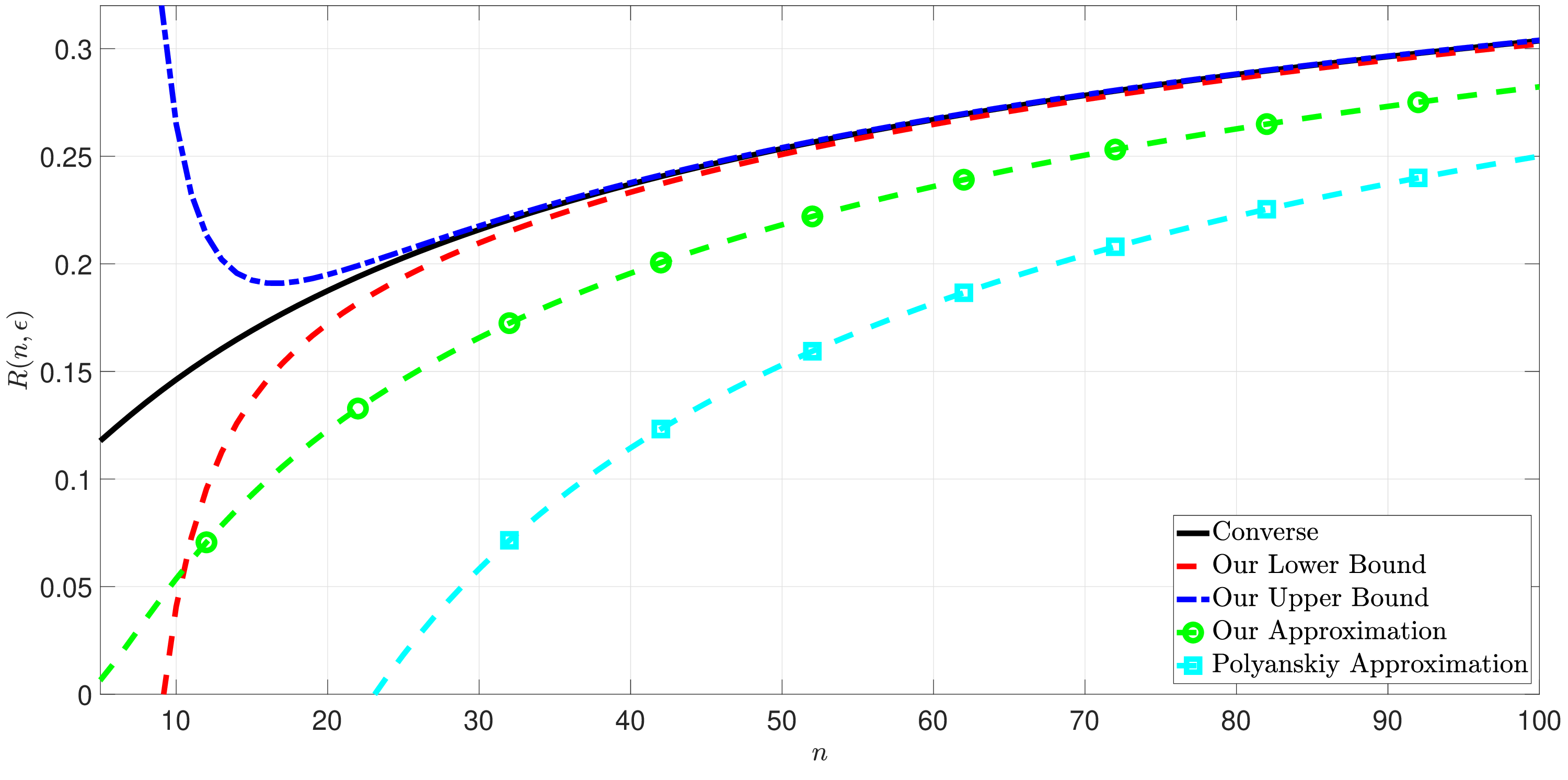}
\caption{Comparison of converse bounds with with bounds and approximatons for $\Omega=1$ and $\epsilon -10^{-3}$}.
\label{fig_l-cap_1}
\end{figure}

\begin{figure}
\includegraphics[width=1\textwidth]{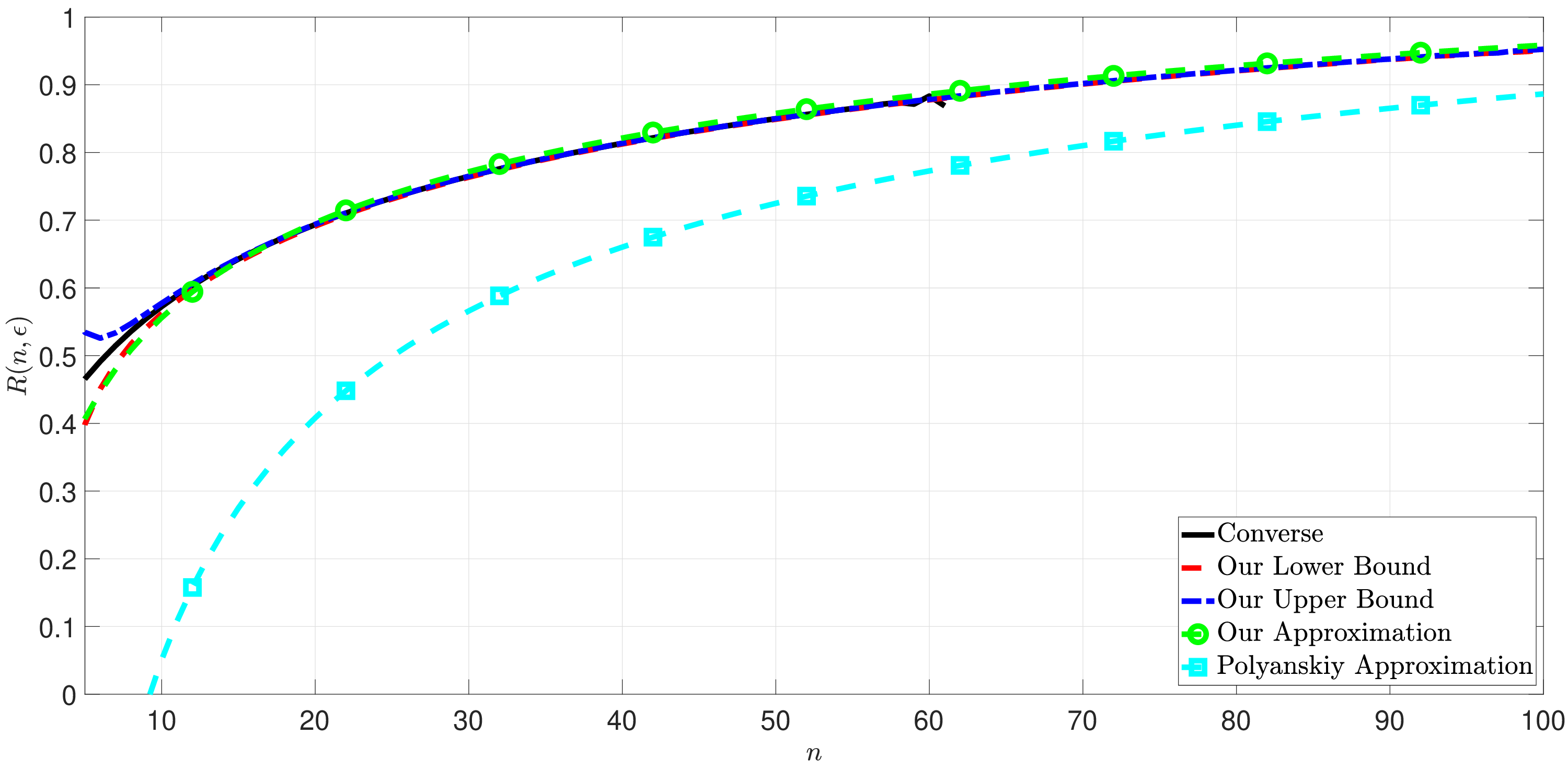}
\caption{Comparison of converse bounds with with bounds and approximatons for $\Omega=5$ and $\epsilon =10^{-5}$}
\label{fig_l-cap_2}
\end{figure}

In this section, we compare the derived upper and lower bounds, as well as the  closed-form approximation, to the converse  bound of the AWGN channel capacity and Polyanskiy-Poor-Verdú closed-form approximation  in the finite-blocklength regime. To this end, the converse  bound is obtained using \eqref{eq:cap_1}, the lower and upper bounds are obtained using \eqref{eq:F_FA_bound1} and  \eqref{eq:F_FA_bound2}, respectively, the closed-form approximation using \eqref{eq:My_rate1} and the Polyanskiy-Poor-Verdú closed-for approximation using \eqref{eq:Polyanskiy_rate}. Two examples are plotted in Figs~\ref{fig_l-cap_1} and \ref{fig_l-cap_2}, one for $\Omega=1$ and $ \epsilon=10^{-3}$ and the other for $\Omega=5$ and  $\epsilon=10^{-5}$, respectively. As  the examples in Figs~\ref{fig_l-cap_1} and \ref{fig_l-cap_2} show, the bounds become tighter as $n$ increases. Moreover, the derived closed-form approximation is in general tighter than the closed-form approximation by Polyanskiy-Poor-Verdú, and its accuracy is better for higher $\Omega$ and  smaller $\epsilon$.

Other possible applications of the derived tail bounds could be to  predict fault-tolerant switching control for discrete linear systems with actuator random failures, see \cite{SHI2024108554}.

\section{Conclusion}\label{sec-conc}

We provided a general method for upper and lower bounding both the left and the right tail probabilities of continuous random variables. The proposed method requires setting a function $g (x)$ with certain conditions, which if satisfied, results in upper and lower bounds that are  functions of  $g (x)$, $g' (x)$, and the PDF of $X$, $f (x)$. We also proposed an iterative method that results in ever tighter upper and lower bounds on the tails, under certain conditions. We established connections between the proposed bounds and Markov's inequality and Chernoff's bound. Finally, we used the proposed bounding method to derive a novel closed-form approximation to the converse bound of the AWGN channel capacity.

\appendices
 
 \section{Proof Of Lemma~\lowercase{\ref{lema_1}}}\label{app_1}

If 
\begin{align}\label{app_1_eq_1}
\frac{1-F (x)}{g (x)} 
\end{align}
is a deceasing function on the interval $I$,  then the following holds
\begin{align}\label{app_1_eq_2}
\diff{ }{ x}\left(\frac{1-F (x)}{g (x)} \right)\leq 0, \; \forall x\in I.
\end{align}
By expanding the derivative in \eqref{app_1_eq_2}, we obtain
\begin{align}\label{app_1_eq_3}
\frac{-f (x)}{g (x)} - \frac{\big(1-F (x)\big) g' (x)}{g^2 (x)}\leq 0 , \; \forall x\in I.
\end{align}
Multiplying both sides of \eqref{app_1_eq_3} by $g^2 (x)$, we obtain
\begin{align}\label{app_1_eq_4}
 -f (x) g (x)  - \big(1-F (x)\big) g' (x) \leq 0 , \; \forall x\in I.
\end{align} 
Now dividing both sides by $g' (x)$, and taking into account that $g' (x)<0$, we obtain
 \begin{align}\label{app_1_eq_5}
 -f (x) \frac{g (x)}{g' (x)}  -  \big(1-F (x)\big)  \geq 0 , \; \forall x\in I.
\end{align} 
Rearranging \eqref{app_1_eq_5}, we obtain the desired result in \eqref{seq_1_eq_3}.
  
On the other hand, if 
\begin{align}\label{app_1_eq_6}
\diff{ }{ x}\left(\frac{1-F (x)}{g (x)} \right)\geq 0, \; \forall x\in I,
\end{align}
then following the same procedure as above leads to the desired result in \eqref{seq_1_eq_6}. This completes the proof.

 \section{Proof Of Theorem~\lowercase{\ref{thm_1}}}\label{app_2}

We start with the bound given by \eqref{seq_1_eq_3} in  Lemma~\ref{lema_1}, which for $I=[\hat x_0,r]$, can be equivalently written as  
\begin{align}\label{app_2_eq_1}
1-F (x)+ f (x)\frac{g (x)}{g' (x)} \leq 0 , \; \forall x>\hat x_0.
\end{align}
As stated in Lemma~\ref{lema_1}, the bound in \eqref{app_2_eq_1} holds if $g (x)$ is continuous, $g (x)>0$, and $g' (x)<0$, $\forall x>\hat x_0$. Now let us define a function $D (x)$ that is equal to the left-hand side of \eqref{app_2_eq_1}, and thereby given by
\begin{align}\label{app_2_eq_2}
D (x)=1-F (x)+ f (x)\frac{g (x)}{g' (x)}.
\end{align}

First note the following obvious property: If a function $G(x)$  is an increasing function for $x_0<x<r$ and if $G(x)$ converges to $\lim\limits_{x\to r} G(x)=0$, then the function $G(x)$ must be a non-positive function for $x_0<x\leq r$, i.e., $G(x)\leq 0$  for $x_0<x\leq r$.

We now use this property for the construction of this proof. Specifically, in the following, we investigate the properties that $g (x)$ must satisfy in order for $D (x)$, given by \eqref{app_2_eq_2}, to satisfy $i)$ $\lim\limits_{x\to r} D (x)=0$ and $ii)$ $D (x)$ to be an increasing function $\forall x>x_0$, since then  $D (x)\leq 0$, $\forall x>x_0$ holds, as per the property described above. On the other hand, when $D (x)\leq 0$, $\forall x>x_0$ holds, then the upper bound in \eqref{seq_1_eq_10} holds, and thereby we have obtained our proof for the   upper bound in \eqref{seq_1_eq_10}.

We start with investigating the conditions of $g (x)$ for which $i)$ $\lim\limits_{x\to r} D (x)=0$  holds. Now, for $D (x)$ given by \eqref{app_2_eq_2}, condition   $\lim\limits_{x\to r} D (x)=0$ is always met since this theorem assumes that \eqref{seq_1_eq_7b} holds. Specifically, we have
\begin{align}\label{app_2_eq_4}
\lim\limits_{x\to r} D (x) &=\lim\limits_{x\to r} \big(1-F (x)\big)+ \lim_{x\to r} f (x)\frac{g (x)}{g' (x)}
\nonumber\\
&= \lim_{x\to r} f (x)\frac{g (x)}{g' (x)}\stackrel{(a)}{=}0,
\end{align}
where $(a)$ follows from the assumption in this theorem that $g (x)$ is such that condition  \eqref{seq_1_eq_7b} is satisfied.

We now continue investigating the conditions of $g (x)$ that make $ii)$ $D (x)$ to be an increasing function $\forall x>x_0$.
For $D (x)$ to be an increasing function $\forall x>x_0$, then the following must hold
\begin{align}\label{app_2_eq_5}
\diff{ }{  x} D (x)\geq 0, \;\forall x>x_0.
\end{align}
Inserting $D (x)$ from \eqref{app_2_eq_2} into \eqref{app_2_eq_5} and carrying out the derivative, we obtain
\begin{align}\label{app_2_eq_6}
&-f (x)+ f' (x)\frac{g (x)}{g' (x)} +f (x)\left(\frac{g' (x)}{g' (x)} - \frac{g (x)g'' (x)}{\big(g' (x)\big)^2}\right)  \geq 0,
\nonumber\\
&  \;\forall x>x_0,
\end{align}
which is equivalent to 
\begin{align}\label{app_2_eq_7}
-f (x)+ f' (x)\frac{g (x)}{g' (x)} +f (x)  -f (x) \frac{g (x)g'' (x)}{\big(g' (x)\big)^2}  \geq 0, \;\forall x>x_0,
\end{align}
which is equivalent to 
\begin{align}\label{app_2_eq_8}
  f' (x)\frac{g (x)}{g' (x)}    -f (x) \frac{g (x)g'' (x)}{\big(g' (x)\big)^2}   \geq 0, \;\forall x>x_0.
\end{align}
Multiplying both sides of \eqref{app_2_eq_8} by $\frac{g (x)}{g' (x)}$, and taking into consideration that $\frac{g (x)}{g' (x)}<0$, $\forall x>x_0$, we obtain the following equivalent inequality
\begin{align}\label{app_2_eq_9}
  f' (x)    -f (x) \frac{g'' (x)}{g' (x)}   \leq 0, \;\forall x>x_0.
\end{align}
Dividing both sides of \eqref{app_2_eq_9} by $f (x)$,  we obtain the following equivalent inequality
\begin{align}\label{app_2_eq_10}
  \frac{f' (x)}{f (x)}    - \frac{g'' (x)}{g' (x)}   \leq 0, \;\forall x>x_0,
\end{align}
which is equivalent to
\begin{align}\label{app_2_eq_11}
 \ln'\big(f (x)\big)    -  \ln'\big(g' (x)\big)  \leq 0, \;\forall x>x_0,
\end{align}
which is equivalent to
\begin{align}\label{app_2_eq_12}
 \ln'\left(\frac{f (x)}{g' (x)}\right)      \leq 0, \;\forall x>x_0,
\end{align}
and also to
\begin{align}\label{app_2_eq_13}
 \ln'\left(\frac{f (x)}{-g' (x)}\right)      \leq 0, \;\forall x>x_0.
\end{align}
The expression \eqref{app_2_eq_13} tells us that $D (x)$ is an increasing function $\forall x>x_0$, if the function
\begin{align}\label{app_2_eq_14}
 \ln\left(\frac{f (x)}{-g' (x)}\right) 
\end{align}
is a decreasing function\footnote{Note that since $g' (x)<0$, the solution of the differential equation in \eqref{app_2_eq_12} or in \eqref{app_2_eq_13} cannot result in $\ln\left(\frac{ f (x)}{g' (x)}\right)$ since this is not a real function, and must result in \eqref{app_2_eq_14}; a result which is obtained by setting the constant of the corresponding differential equation such that the solution is a real function.},  $\forall x>x_0$. 
Now since $\ln(\cdot)$ is a one-to-one function, \eqref{app_2_eq_14} is a  decreasing function when
\begin{align}\label{app_2_eq_15}
 \frac{f (x)}{-g' (x)} 
\end{align}
is a decreasing function $\forall x>x_0$. We now simplify the condition that \eqref{app_2_eq_15} is decreasing function $\forall x>x_0$ as follows. The function in \eqref{app_2_eq_15} is a deceasing function $\forall x>x_0$ if
\begin{align}\label{app_2_eq_15a}
 \diff{ }{  x }\left(\frac{f (x)}{-g' (x)} \right)\leq 0, \;\forall x>x_0.
\end{align}
On the other hand, the left-hand side of \eqref{app_2_eq_15a} can be written equivalently as
\begin{align}\label{app_2_eq_15b}
  \diff{ }{  x } \left(\frac{f (x)}{-g' (x)} \right)  &=   \diff{ }{  x }\left( \frac{f (x) g (x)}{-g' (x)g (x)}  \right)
 =   \diff{ }{  x } \left( \frac{P (x)}{g (x)}  \right)
 \nonumber\\
 & 
 = \frac{P' (x)}{g (x)}-\frac{P (x)g' (x)}{g^2 (x)} 
  \nonumber\\
 &  = \frac{P' (x)}{g (x)}-\frac{P (x)g' (x)f (x)}{g^2 (x)f (x)}
  \nonumber\\
 & 
= \frac{P' (x)}{g (x)}+\frac{P (x)f (x)}{P (x)g (x)}
= \frac{P' (x)}{g (x)}+\frac{f (x)}{g (x)}.
\end{align}
Inserting \eqref{app_2_eq_15b} into \eqref{app_2_eq_15a} and multiplying both sides of the inequality by $g (x)$, we obtain
\begin{align}\label{app_2_eq_15c}
  P' (x) + f (x) \leq 0,\; \forall x>x_0.
\end{align}
Therefore, if  \eqref{app_2_eq_15c} holds $\forall x>x_0$ and since \eqref{app_2_eq_4} always holds, then 
 \begin{align}\label{app_2_eq_16}
 D (x)\leq 0, \; \forall x>x_0.
\end{align}
Inserting \eqref{app_2_eq_2} into \eqref{app_2_eq_16}, we obtain \eqref{seq_1_eq_10}, which is the first bound we aimed to prove.

For the second bound in this theorem, we   follow the same method as above.   
To this end, we use the following obvious property: If a function $G(x)$  is a decreasing function for $x_0<x<r$ and if $G(x)$ converges to $\lim\limits_{x\to r} G(x)=0$, then the function $G(x)$ must be a non-negative function for $x_0<x\leq r$, i.e., $G(x)\geq 0$ for $x_0<x\leq r$.  
 Thereby, it is straightforward to prove that   $D (x)\geq 0$ if
\begin{align}\label{app_2_eq_17}
 \frac{f (x)}{-g' (x)} 
\end{align}
is an increasing function $\forall x>x_0$. Taking into account the   proof of the first bound, the proof of the second bound is omitted due to its redundancy. This concludes the proof of this theorem.

 \section{Proof Of Lemma~\lowercase{\ref{lema_2}}}\label{app_2-3}  
Since  $\lim\limits_{x\to r} g (x)=0$ holds, we have the following limit for $P (x)$
 \begin{align}\label{app_2-3_eq_1}
 \lim_{x\to r}P (x) &= \lim_{x\to r} - f (x)\frac{g (x)}{g' (x)} = \lim_{x\to r}  g (x) \lim_{x\to r}\frac{-f (x)}{g' (x)}\nonumber\\
& 
\stackrel{(a)}{=}\lim_{x\to r}  g (x) \lim_{x\to r}\frac{1-F (x)}{g (x)}\nonumber\\
&=\lim_{x\to r}  g (x)  \frac{1-F (x)}{g (x)}=\lim_{x\to r}   \big(1-F (x)\big)=0,
\end{align}
where $(a)$ follows from (reverse) l'Hopital's rule, which is valid since $1-F (x)$ and $g (x)$ both  tend to $0$ as $x\to r$.

 \section{Proof Of Lemma~\lowercase{\ref{lema_2a}}}\label{app_2-3a}  
 
 Since we assume that \eqref{seq_2_eq_1} is valid upper bound or a valid lower bound, $\forall x>x_0$, as per Theorem~\ref{thm_1}, then $\lim\limits_{x\to r} P_0(x)=0$ must hold. Now, if $\lim\limits_{x\to r} P_0(x)=0$ holds, then the next iteration, $P_1(x)$, obtained by setting $i=0$ in \eqref{seq_2_eq_2}, must satisfy 
 \begin{align}\label{seq_2_eq_2-2}
\lim\limits_{x\to r} P_{1}(x)=0,
\end{align}
 according to Lemma~\ref{lema_2}.
Now, since  \eqref{seq_2_eq_2-2} holds, then $P_2(x)$, obtained by setting $i=1$ in \eqref{seq_2_eq_2} must also satisfy 
 \begin{align}\label{seq_2_eq_2-3}
\lim\limits_{x\to r} P_{2}(x)=0,
\end{align}
 according to Lemma~\ref{lema_2}.
These true statements  can be extended to any $i+1$, thereby proving that \eqref{seq_1_eq_13-n2-1} holds.

 \section{Proof Of Lemma~\lowercase{\ref{lema_3}}}\label{app_3} 
 
 Let $P_i$ be an upper bound on $1-F (x)$, $\forall x>x_0$, as per Theorem~\ref{thm_1}. Then, according to \eqref{seq_1_eq_9}, the following holds
 \begin{align}\label{app_3_eq_1}
P'_i(x)\leq-f (x),\; \forall x>x_0,
\end{align}
which is equivalent to
 \begin{align}\label{app_3_eq_2}
1\geq-f (x)\frac{1}{P'_i(x)},\; \forall x>x_i,
\end{align}
and is obtained by  dividing both sides of \eqref{app_3_eq_2} by $P'_i(x)<0$, which is a valid division only for $x>\hat x$, due to the assumption $P'_i(x)<0$, $x>\hat x$.
Now, \eqref{app_3_eq_2}
is equivalent to 
 \begin{align}\label{app_3_eq_3}
1\geq-f (x)\frac{P_i(x)}{P'_i(x)}\frac{1}{P_i(x)},\; \forall x>x_i,
\end{align}
which is equivalent to
 \begin{align}\label{app_3_eq_4}
1\geq \frac{P_{i+1}(x)}{P_i(x)},\; \forall x>x_i,
\end{align}
which is equivalent to
 \begin{align}\label{app_3_eq_5}
P_i(x)\geq  P_{i+1}(x) ,\; \forall x>x_i,
\end{align}
which is the result in \eqref{seq_2_eq_3} that we aimed to prove.

We now prove \eqref{seq_2_eq_4}.
Let $P_i$ be a lower bound on $1-F (x)$, $\forall x>x_0$, as per Theorem~\ref{thm_1}. Then, according to \eqref{seq_1_eq_12}, the following holds
 \begin{align}\label{app_3_eq_6}
P'_i(x)\geq-f (x),\; \forall x>x_0,
\end{align}
which is equivalent to
 \begin{align}\label{app_3_eq_7}
1\leq-f (x)\frac{1}{P'_i(x)},\; \forall x>x_i,
\end{align}
and is obtained by  dividing both sides of \eqref{app_3_eq_7} by $P'_i(x)<0$, which is a valid division only for $x>\hat x$. Now, \eqref{app_3_eq_7}
is equivalent to 
 \begin{align}\label{app_3_eq_8}
1\leq-f (x)\frac{P_i(x)}{P'_i(x)}\frac{1}{P_i(x)},\; \forall x>x_i,
\end{align}
which is equivalent to
 \begin{align}\label{app_3_eq_9}
1\leq \frac{P_{i+1}(x)}{P_i(x)},\; \forall x>x_i,
\end{align}
which is equivalent to
 \begin{align}\label{app_3_eq_10}
P_i(x)\leq  P_{i+1}(x) ,\; \forall x>x_i,
\end{align}
which is the result in \eqref{seq_2_eq_4} that we aimed to prove.

 \section{Proof Of Lemma~\lowercase{\ref{lema_4}}}\label{app_4}

We start with the first part of this lemma when $ P_{i}(x)$ and $ P_{i+1}(x)$   are both valid upper and lower bounds on $1-F (x)$, $\forall x>x_{i+1}$, respectively, as per Theorem~\ref{thm_1}. 

First, note that from \eqref{seq_2_eq_4a}, that the following holds
  \begin{align}\label{app_4_eq_1}
\frac{P_{L,i}(x)+ P_i(x)}{2} =1-F (x), \;\forall x>x_{i+1}.
\end{align}
Then, if $P_{L,i}(x)\leq P_{i+1}(x)$ holds, from \eqref{app_4_eq_1}, we obtain that the following  must  hold
  \begin{align}\label{app_4_eq_2}
\frac{P_{i+1}(x)+ P_i(x)}{2} \geq 1-F (x), \;\forall x>x_{i+1},
\end{align}
which means that $\frac{P_{i+1}(x)+ P_i(x)}{2} $ is an upper bound on $1-F (x)$, $\forall x>x_{i+1}$. Now, according to Theorem~\ref{thm_1}, in order for   $\frac{P_{i+1}(x)+ P_i(x)}{2} $ to  be an upper bound on $1-F (x)$, $\forall x>x_{i+1}$, i.e., \eqref{app_4_eq_2} to hold, the following must hold
  \begin{align}\label{app_4_eq_3}
\frac{P'_{i+1}(x)+ P'_i(x)}{2} +f (x) \leq 0, \;\forall x>x_{i+1},
\end{align}
which is equivalent to
  \begin{align}\label{app_4_eq_4}
P'_{i+1}(x)+ P'_i(x) + 2f (x) \leq 0,\; \forall x>x_{i+1}.
\end{align}
Going now in reverse order, it is straightforward to obtain that if \eqref{app_4_eq_4} holds, then $P_{L,i}(x)\leq P_{i+1}(x)$ also holds. This completes the first part of this lemma.

We now prove the second part of this lemma when $ P_{i}(x)$ and $ P_{i+1}(x)$   are both valid lower and upper bounds on $1-F (x)$, $\forall x>x_{i+1}$, respectively, as per Theorem~\ref{thm_1}. 

First, note that from \eqref{seq_2_eq_4aa}, that the following holds
  \begin{align}\label{app_4_eq_5}
\frac{P_{U,i}(x)+ P_i(x)}{2} =1-F (x),\; \forall x>x_{i+1}.
\end{align}
Then, if $P_{U,i}(x)\geq P_{i+1}(x)$ holds, from \eqref{app_4_eq_5}, we obtain that the following must  hold
  \begin{align}\label{app_4_eq_6}
\frac{P_{i+1}(x)+ P_i(x)}{2} \leq 1-F (x), \; \forall x>x_{i+1},
\end{align}
which means that $\frac{P_{i+1}(x)+ P_i(x)}{2} $ is a lower bound  on $1-F (x)$, $\forall x>x_{i+1}$. Now, according to Theorem~\ref{thm_1}, in order for   $\frac{P_{i+1}(x)+ P_i(x)}{2} $ to  be an lower bound  on $1-F (x)$, $\forall x>x_{i+1}$, i.e., \eqref{app_4_eq_6} to hold, the following must hold
  \begin{align}\label{app_4_eq_7}
\frac{P'_{i+1}(x)+ P'_i(x)}{2} +f (x) \geq 0, \;\forall x>x_{i+1},
\end{align}
which is equivalent to
  \begin{align}\label{app_4_eq_8}
P'_{i+1}(x)+ P'_i(x) + 2f (x) \geq 0,\; \forall x>x_{i+1}.
\end{align}
Going now in reverse order, it is straightforward to obtain that if \eqref{app_4_eq_8} holds, then $P_{U,i}(x)\geq P_{i+1}(x)$ also holds.  This completes the second part of this lemma.

Finally, note that we have used Theorem~\ref{thm_1} for determining whether $\frac{P_{i+1}(x)+ P_i(x)}{2}$ is an upper bound or a lower bound, as per Theorem~\ref{thm_1}, and this  is valid since
\begin{align}\label{app_4_eq_9}
\lim_{x\to r} \frac{P_{i+1}(x)+ P_i(x)}{2}=0
\end{align}
holds   due to Lemma~\ref{lema_2a}. This completes the proof  of this lemma.

 \section{Proof Of Theorem~\lowercase{\ref{thm_2}}}\label{app_5}

We aim to prove  
\begin{align}\label{app_5_eq_0}
 F (x) \leq   f (x)\frac{g (x)}{g' (x)} , \; \forall x<x_0,
\end{align}
which is equivalent to 
\begin{align}\label{app_5_eq_1}
 F (x) -  f (x)\frac{g (x)}{g' (x)} \leq 0 , \; \forall x<x_0.
\end{align}
 Now let us again define a function $D (x)$ that is equal to the left-hand side of \eqref{app_5_eq_1}, and thereby given by
\begin{align}\label{app_5_eq_2}
D (x)= F (x)- f (x)\frac{g (x)}{g' (x)}.
\end{align}

Note now the following obvious property: If a function $G(x)$  is a  decreasing function for $l<x<x_0$ and if $G(x)$ converges to $\lim\limits_{x\to l} G(x)=0$, then the function $G(x)$ must be a non-positive function for $l<x\leq x_0$, i.e., $G(x)\leq 0$ for $l<x\leq x_0$.

We now use this property for the construction of this proof. Specifically, in the following, we investigate the properties that $g (x)$ must satisfy in order for $D (x)$, given by \eqref{app_5_eq_2}, to satisfy $i)$ $\lim\limits_{x\to l} D (x)=0$ and $ii)$ $D (x)$ to be a decreasing function $\forall x<x_0$, since then  $D (x)\leq 0$, $\forall x<x_0$ holds, according to property described above. On the other hand, when $D (x)\leq 0$, $\forall x<x_0$, holds, then the upper bound in \eqref{seq_5_eq_10} holds, and thereby we have obtained our proof for the  upper bound in \eqref{seq_5_eq_10}.

We start with investigating the conditions of $g (x)$ for which $i)$ $\lim\limits_{x\to l} D (x)=0$  holds. Now, for $D (x)$ given by \eqref{app_5_eq_2}, condition   $\lim\limits_{x\to l} D (x)=0$ is always met since this theorem assumes that \eqref{seq_5_eq_7b} holds. Specifically, we have
\begin{align}\label{app_5_eq_4}
\lim\limits_{x\to l} D (x) &=\lim\limits_{x\to l}  F (x) - \lim_{x\to l} f (x)\frac{g (x)}{g' (x)}= -\lim_{x\to l} f (x)\frac{g (x)}{g' (x)}\stackrel{(a)}{=}0,
\end{align}
where $(a)$ follows from the assumption in this theorem that $g (x)$ is such that condition  \eqref{seq_5_eq_7b} is satisfied.

We now continue investigating the conditions of $g (x)$ that make $ii)$ $D (x)$ to be a  decreasing function $\forall x<x_0$.
For $D (x)$ to be a  decreasing function $\forall x<x_0$,   the following must hold
\begin{align}\label{app_5_eq_5}
\diff{ }{  x} D (x)\leq 0, \;\forall x<x_0.
\end{align}
Inserting $D (x)$ from \eqref{app_5_eq_2} into \eqref{app_5_eq_5} and carrying out the derivative, we obtain
\begin{align}\label{app_5_eq_6}
 f (x)- f' (x)\frac{g (x)}{g' (x)} -f (x)\left(\frac{g' (x)}{g' (x)} - \frac{g (x)g'' (x)}{\big(g' (x)\big)^2}\right)  \leq 0,  \;\forall x<x_0,
\end{align}
which is equivalent to 
\begin{align}\label{app_5_eq_7}
 f (x)- f' (x)\frac{g (x)}{g' (x)} -f (x)  +f (x) \frac{g (x)g'' (x)}{\big(g' (x)\big)^2}  \leq 0, \;\forall x<x_0,
\end{align}
which is equivalent to 
\begin{align}\label{app_5_eq_8}
 - f' (x)\frac{g (x)}{g' (x)}    +f (x) \frac{g (x)g'' (x)}{\big(g' (x)\big)^2}   \leq 0, \;\forall x<x_0.
\end{align}
Multiplying both sides of \eqref{app_5_eq_8} by $-\frac{g (x)}{g' (x)}$, and taking into consideration that $-\frac{g (x)}{g' (x)}<0$, $\forall x<x_0$, we obtain the following equivalent inequality
\begin{align}\label{app_5_eq_9}
  f' (x)    - f (x) \frac{g'' (x)}{g' (x)}   \geq 0, \;\forall x<x_0.
\end{align}
Dividing both sides of \eqref{app_5_eq_9} by $f (x)$,  we obtain the following equivalent inequality
\begin{align}\label{app_5_eq_10}
   \frac{f' (x)}{f (x)}    - \frac{g'' (x)}{g' (x)}   \geq 0, \;\forall x<x_0,
\end{align}
which is equivalent to
\begin{align}\label{app_5_eq_11}
 \ln'\big(f (x)\big)    -  \ln'\big(g' (x)\big)  \geq 0, \;\forall x<x_0,
\end{align}
which is equivalent to
\begin{align}\label{app_5_eq_12}
 \ln'\left(\frac{f (x)}{g' (x)}\right)      \geq 0, \;\forall x<x_0,
\end{align}
and also to
\begin{align}\label{app_5_eq_13}
 \ln'\left(\frac{f (x)}{-g' (x)}\right)      \geq 0, \;\forall x<x_0.
\end{align}
The expression \eqref{app_5_eq_12} tells us that $D (x)$ is a  decreasing function $\forall x<x_0$ if the function
\begin{align}\label{app_5_eq_14}
 \ln\left(\frac{f (x)}{ g' (x)}\right) 
\end{align}
is an increasing function\footnote{Note that since $g' (x)>0$, the solution of the differential equation in \eqref{app_5_eq_13} or in \eqref{app_5_eq_12} cannot result in $\ln\left(\frac{ f (x)}{-g' (x)}\right)$ since this is not a real function, and must result in \eqref{app_5_eq_14}; a result which is obtained by setting the constant of the corresponding differential equation such that the solution is a real function.},  $\forall x<x_0$. 
Now since $\ln(\cdot)$ is a one-to-one function, \eqref{app_5_eq_14} is an increasing function when
\begin{align}\label{app_5_eq_15}
 \frac{f (x)}{g' (x)} 
\end{align}
is an increasing function $\forall x<x_0$. We now simplify the condition that \eqref{app_5_eq_15} is an increasing function $\forall x<x_0$ as follows. The function in \eqref{app_5_eq_15} is an increasing function $\forall x<x_0$ if
\begin{align}\label{app_5_eq_15a}
 \diff{ }{  x }\left(\frac{f (x)}{g' (x)} \right)\geq 0, \;\forall x<x_0.
\end{align}
On the other hand, the left-hand side of \eqref{app_5_eq_15a} can be written equivalently as
\begin{align}\label{app_5_eq_15b}
 \diff{ }{  x }\left( \frac{f (x)}{ g' (x)}\right) &=   \diff{ }{  x }\left(\frac{f (x) g (x)}{ g' (x)g (x)}\right)  =   \diff{ }{  x }\left(\frac{P (x)}{g (x)}\right)
 \nonumber\\
 &  
  = \frac{P' (x)}{g (x)}-\frac{P (x)g' (x)}{g^2 (x)}\nonumber\\
& = \frac{P' (x)}{g (x)}-\frac{P (x)g' (x)f (x)}{g^2 (x)f (x)}
 \nonumber\\
 & 
= \frac{P' (x)}{g (x)}-\frac{P (x)f (x)}{P (x)g (x)}
= \frac{P' (x)}{g (x)}-\frac{f (x)}{g (x)}.
\end{align}
Inserting \eqref{app_5_eq_15b} into \eqref{app_5_eq_15a} and multiplying both sides of the inequality by $g (x)$, we obtain
\begin{align}\label{app_5_eq_15c}
  P' (x) - f (x) \geq 0,\;\forall x<x_0.
\end{align}
Therefore, if  \eqref{app_5_eq_15c} holds $\forall x<x_0$ and since \eqref{app_5_eq_4} always holds, then 
 \begin{align}\label{app_5_eq_16}
 D (x)\leq 0, \; \forall x<x_0.
\end{align}
Inserting \eqref{app_5_eq_2} into \eqref{app_5_eq_16}, we obtain \eqref{seq_5_eq_10}, which is the first bound we aimed to prove.

For the second bound in this Theorem, we   follow the same method as above. This this end, we use the following obvious property: If a function $G(x)$  is an increasing function for $l<x<x_0$ and if $G(x)$ converges to $\lim\limits_{x\to l} G(x)=0$, then the function $G(x)$ must be a non-negative function for $l<x\leq x_0$,  , i.e., $G(x)\geq 0$ for $l<x\leq x_0$. Thereby, it is straightforward to prove that   $D (x)\geq 0$ if
\begin{align}\label{app_5_eq_17}
 \frac{f (x)}{g' (x)} 
\end{align}
is a decreasing function $\forall x<x_0$. Taking into account the   proof of the first bound, the proof of the second bound is omitted due to its redundancy. This concludes the proof of this theorem.

 \section{Proof Of Lemma~\lowercase{\ref{lema_6}}}\label{app_7} 

Since  $\lim\limits_{x\to l} g (x)=0$ holds, we have the following limit for $P (x)$
 \begin{align}\label{app_7_eq_1}
 \lim_{x\to l}P (x) &= \lim_{x\to l}   f (x)\frac{g (x)}{g' (x)} = \lim_{x\to l}  g (x) \lim_{x\to l}\frac{ f (x)}{g' (x)}
  \nonumber\\
 & 
\stackrel{(a)}{=}\lim_{x\to l}  g (x) \lim_{x\to l}\frac{F (x)}{g (x)}\nonumber\\
&=\lim_{x\to l}  g (x)  \frac{F (x)}{g (x)}=\lim_{x\to l}  F (x)=0,
\end{align}
where $(a)$ follows from (reverse) l'Hopital's rule, which is valid since $F (x)$ and $g (x)$ both  tend to $0$ as $x\to l$.

 \section{Proof Of Lemma~\lowercase{\ref{lema_6a}}}\label{app_8}  
 
 Since we assume that \eqref{seq_6_eq_1} is valid upper bound or a valid lower bound, $\forall x<x_0$, as per Theorem~\ref{thm_2}, then $\lim\limits_{x\to l} P_0(x)=0$ must hold. Now, if $\lim\limits_{x\to l} P_0(x)=0$ holds, then the next iteration, $P_1(x)$, obtained by setting $i=0$ in \eqref{seq_6_eq_2}, must satisfy 
 \begin{align}\label{seq_6_eq_2-2}
\lim\limits_{x\to l} P_{1}(x)=0,
\end{align}
 according to Lemma~\ref{lema_6}.
Now, since  \eqref{seq_6_eq_2-2} holds, then $P_2(x)$, obtained by setting $i=1$ in \eqref{seq_6_eq_2} must also satisfy 
 \begin{align}\label{seq_6_eq_2-3}
\lim\limits_{x\to l} P_{2}(x)=0,
\end{align}
 according to Lemma~\ref{lema_6}.
These true statements  can be extended to any $i+1$, thereby proving that \eqref{seq_5_eq_13-n2-1} holds.

 \section{Proof Of Lemma~\lowercase{\ref{lema_7}}}\label{app_9} 
 
 Let $P_i$ be an upper bound on $F (x)$, $\forall x<x_0$, as per Theorem~\ref{thm_2}. Then, according to \eqref{seq_5_eq_9}, the following holds
 \begin{align}\label{app_9_eq_1}
P'_i(x)\geq f (x),\; \forall x<x_0,
\end{align}
which is equivalent to
 \begin{align}\label{app_9_eq_2}
1\geq f (x)\frac{1}{P'_i(x)},\; \forall x<x_i,
\end{align}
and is obtained by  dividing both sides of \eqref{app_9_eq_2} by $P'_i(x)>0$, which is a valid division only for $x<\hat x$, due to the assumption $P'_i(x)>0$, $x<\hat x$.
Now, \eqref{app_9_eq_2}
is equivalent to 
 \begin{align}\label{app_9_eq_3}
1\geq f (x)\frac{P_i(x)}{P'_i(x)}\frac{1}{P_i(x)},\; \forall x<x_i,
\end{align}
which is equivalent to
 \begin{align}\label{app_9_eq_4}
1\geq \frac{P_{i+1}(x)}{P_i(x)},\; \forall x<x_i,
\end{align}
which is equivalent to
 \begin{align}\label{app_9_eq_5}
P_i(x)\geq  P_{i+1}(x) ,\; \forall x<x_i,
\end{align}
which is the result in \eqref{seq_6_eq_3} that we aimed to prove.

We now prove \eqref{seq_6_eq_4}.
Let $P_i$ be a lower bound on $F (x)$, $\forall x<x_0$, as per Theorem~\ref{thm_2}. Then, according to \eqref{seq_5_eq_12}, the following holds
 \begin{align}\label{app_9_eq_6}
P'_i(x)\leq f (x),\; \forall x<x_0,
\end{align}
which is equivalent to
 \begin{align}\label{app_9_eq_7}
1\leq f (x)\frac{1}{P'_i(x)},\; \forall x<x_i,
\end{align}
and is obtained by  dividing both sides of \eqref{app_9_eq_7} by $P'_i(x)>0$, which is a valid division only for $x<\hat x$. Now, \eqref{app_9_eq_7}
is equivalent to 
 \begin{align}\label{app_9_eq_8}
1\leq f (x)\frac{P_i(x)}{P'_i(x)}\frac{1}{P_i(x)},\; \forall x<x_i,
\end{align}
which is equivalent to
 \begin{align}\label{app_9_eq_9}
1\leq \frac{P_{i+1}(x)}{P_i(x)},\; \forall x<x_i,
\end{align}
which is equivalent to
 \begin{align}\label{app_9_eq_10}
P_i(x)\leq  P_{i+1}(x) ,\; \forall x<x_i,
\end{align}
which is the result in \eqref{seq_6_eq_4} that we aimed to prove.

 \section{Proof Of Lemma~\lowercase{\ref{lema_8}}}\label{app_10}

We start with the first part of this lemma when $ P_{i}(x)$ and $ P_{i+1}(x)$   are both valid upper and lower bounds on $F (x)$, $\forall x<x_{i+1}$, respectively, as per Theorem~\ref{thm_2}. 

First, note that from \eqref{seq_6_eq_4a}, that the following holds
  \begin{align}\label{app_10_eq_1}
\frac{P_{L,i}(x)+ P_i(x)}{2} =F (x), \;\forall x<x_{i+1}.
\end{align}
Then, if $P_{L,i}(x)\leq P_{i+1}(x)$ holds, from \eqref{app_10_eq_1}, we obtain that the following  must  hold
  \begin{align}\label{app_10_eq_2}
\frac{P_{i+1}(x)+ P_i(x)}{2} \geq F (x), \;\forall x<x_{i+1},
\end{align}
which means that $\frac{P_{i+1}(x)+ P_i(x)}{2} $ is an upper bound on $F (x)$, $\forall x<x_{i+1}$. Now, according to Theorem~\ref{thm_2}, in order for   $\frac{P_{i+1}(x)+ P_i(x)}{2} $ to  be an upper bound on $F (x)$, $\forall x<x_{i+1}$, i.e., \eqref{app_10_eq_2} to hold, the following must hold
  \begin{align}\label{app_10_eq_3}
\frac{P'_{i+1}(x)+ P'_i(x)}{2} -f (x) \geq 0, \;\forall x<x_{i+1},
\end{align}
which is equivalent to
  \begin{align}\label{app_10_eq_4}
P'_{i+1}(x)+ P'_i(x) - 2f (x) \geq 0,\; \forall x<x_{i+1}.
\end{align}
Going now in reverse order, it is straightforward to obtain that if \eqref{app_10_eq_4} holds, then $P_{L,i}(x)\leq P_{i+1}(x)$ also holds. This completes the first part of this lemma.

We now prove the second part of this lemma when $ P_{i}(x)$ and $ P_{i+1}(x)$   are both valid lower and upper bounds on $F (x)$, $\forall x<x_{i+1}$, respectively, as per Theorem~\ref{thm_2}. 

First, note that from \eqref{seq_6_eq_4aa}, that the following holds
  \begin{align}\label{app_10_eq_5}
\frac{P_{U,i}(x)+ P_i(x)}{2} =F (x),\; \forall x<x_{i+1}.
\end{align}
Then, if $P_{U,i}(x)\geq P_{i+1}(x)$ holds, from \eqref{app_10_eq_5}, we obtain that the following must  hold
  \begin{align}\label{app_10_eq_6}
\frac{P_{i+1}(x)+ P_i(x)}{2} \leq F (x), \; \forall x<x_{i+1},
\end{align}
which means that $\frac{P_{i+1}(x)+ P_i(x)}{2} $ is a lower bound  on $F (x)$, $\forall x<x_{i+1}$. Now, according to Theorem~\ref{thm_2}, in order for   $\frac{P_{i+1}(x)+ P_i(x)}{2} $ to  be an lower bound  on $F (x)$, $\forall x<x_{i+1}$, i.e., \eqref{app_10_eq_6} to hold, the following must hold
  \begin{align}\label{app_10_eq_7}
\frac{P'_{i+1}(x)+ P'_i(x)}{2} -f (x) \leq 0, \;\forall x<x_{i+1},
\end{align}
which is equivalent to
  \begin{align}\label{app_10_eq_8}
P'_{i+1}(x)+ P'_i(x) - 2f (x) \leq 0,\; \forall x<x_{i+1}.
\end{align}
Going now in reverse order, it is straightforward to obtain that if \eqref{app_10_eq_8} holds, then $P_{U,i}(x)\geq P_{i+1}(x)$ also holds.  This completes the second part of this lemma.

Finally, note that we have used Theorem~\ref{thm_2} for determining whether $\frac{P_{i+1}(x)+ P_i(x)}{2}$ is an upper bound or a lower bound, as per Theorem~\ref{thm_2}, and this  is valid since
\begin{align}\label{app_10_eq_9}
\lim_{x\to l} \frac{P_{i+1}(x)+ P_i(x)}{2}=0
\end{align}
holds   due to Lemma~\ref{lema_6a}. This completes the proof  of this lemma.

\section{Derivation Of The Asymptotic Expression of $\lambda$}\label{sec-app-lam}

We need to solve, for $\lambda>0$,
\begin{align}
&-\,\frac{e^{-\frac{n(1+\lambda\Omega)}{2\Omega}}\;
n\sqrt{\lambda/\Omega}\;(\lambda\Omega)^{n/4}\;
I_{\frac{n}{2}-1}\!\bigl(n\sqrt{\lambda/\Omega}\,\bigr)^2}{\;
(-2+n-n\lambda)\, I_{\frac{n}{2}-1}\!\bigl(n\sqrt{\lambda/\Omega}\,\bigr)
+ n\sqrt{\lambda/\Omega}\, I_{\frac{n}{2}}\!\bigl(n\sqrt{\lambda/\Omega}\,\bigr)}\nonumber\\
&=\epsilon,
\label{eq:main}
\end{align}
where $I_\nu$ is the modified Bessel function of the first kind when $n$ is large.
Introduce the standard large-order scaling
\[
\nu=\frac{n}{2},\qquad x=n\sqrt{\frac{\lambda}{\Omega}},\qquad
z=\frac{x}{\nu}=2\sqrt{\frac{\lambda}{\Omega}},
\]
and write the uniform Debye (large $\nu$) approximation for fixed $z>0$ as
\begin{equation}
I_{\nu}(\nu z)\sim \frac{1}{\sqrt{2\pi \nu}}\,
\frac{e^{\nu \eta(z)}}{(1+z^2)^{1/4}},
\label{eq:debye}
\end{equation}
where
\begin{equation}
\eta(z)=\sqrt{1+z^2}+\ln\!\frac{z}{1+\sqrt{1+z^2}}.
\label{eq:debye1}
\end{equation}
A well known standard corollary is the large-order ratio
\begin{equation}
\frac{I_\nu(\nu z)}{I_{\nu-1}(\nu z)} \;\sim\;
\frac{z}{\,1+\sqrt{1+z^2}\,}\qquad(\nu\to\infty).
\label{eq:ratio}
\end{equation}

We now
rewrite \eqref{eq:main} using $x=\nu z$ and factor one copy of $I_{\nu-1}$ from the denominator
\[
\text{LHS}=
-\,\frac{\sqrt{\lambda/\Omega}\;(\lambda\Omega)^{n/4}\,e^{-\frac{n(1+\lambda\Omega)}{2\Omega}}\;
I_{\nu-1}(x)}{(1-\lambda)+\sqrt{\lambda/\Omega}\,\dfrac{I_\nu(x)}{I_{\nu-1}(x)}}
\cdot \frac{1}{n}\,.
\]
Applying \eqref{eq:debye} to $I_{\nu-1}(x)$ but \emph{with the common exponential $e^{\nu\eta(z)}$} factored for both $I_{\nu-1}$ and $I_\nu$, one finds
\begin{align}
\text{LHS}\sim
&-\,\frac{\sqrt{\lambda/\Omega}}{(1-\lambda)+\sqrt{\lambda/\Omega}\,R(z)}
\cdot \frac{1}{\sqrt{\pi n}}\;\frac{e^{n\Phi(\lambda)}}{(1+z^2)^{1/4}}\nonumber\\
&\times
\frac{1+\sqrt{1+z^2}}{z},
\label{eq:LHS-structure}
\end{align}
where
\begin{equation}\label{eq:R(z)-def}
R(z):=\frac{I_\nu(\nu z)}{I_{\nu-1}(\nu z)}\sim \frac{z}{1+\sqrt{1+z^2}}\quad\text{by \eqref{eq:ratio},}
\end{equation}
\begin{equation}
\Phi(\lambda):=-\frac{1+\lambda\Omega}{2\Omega}+\frac{1}{4}\ln(\lambda\Omega)+\frac{1}{2}\eta(z),
\label{eq:Phi-def}
\end{equation}
and $z=2\sqrt{\lambda/\Omega}$.

Since $0<\epsilon<1$ is $O(1)$, the exponential factor $e^{n\Phi(\lambda)}$ must not blow up or vanish as $n\to\infty$. Thus the leading-order $\lambda=\lambda_0$ is determined by
\begin{equation}
\Phi(\lambda_0)=0.
\label{eq:Phi0}
\end{equation}
Using \eqref{eq:Phi-def} and $z=2\sqrt{\lambda/\Omega}$, \eqref{eq:Phi0} is equivalent to
\begin{equation}
\sqrt{1+\frac{4\lambda}{\Omega}}-\lambda-\frac{1}{\Omega}
+\ln\!\Biggl(\frac{2\lambda}{1+\sqrt{1+\frac{4\lambda}{\Omega}}}\Biggr)=0.
\label{eq:scalar}
\end{equation}
which leads to
\begin{equation}
\lambda_0=1+\frac{1}{\Omega}
\label{eq:lambda0}
\end{equation}

We will now expand both the exponential and the (algebraic) denominator in \eqref{eq:LHS-structure}.

\smallskip
\emph{(i) Exponential.)} Differentiate \eqref{eq:Phi-def}. Using $\eta'(z)=\sqrt{1+z^2}/z$ and
$z=2\sqrt{\lambda/\Omega}$, one finds
\[
\Phi'(\lambda)=-\frac{1}{2}+\frac{1}{4\lambda}+\frac{1}{4\lambda}\sqrt{1+\frac{4\lambda}{\Omega}},
\]
so that $\Phi'(\lambda_0)=0$. A further differentiation gives
\[
\Phi''(\lambda_0)=-\frac{\Omega}{2(\Omega+2)}<0.
\]
Hence, with the scaled displacement
\begin{equation}
u:=\sqrt{n}\,(\lambda-\lambda_0),
\label{eq:u-def}
\end{equation}
we have the Gaussian expansion
\begin{equation}
e^{n\Phi(\lambda)}=\exp\!\left(-\,\frac{\Omega}{4(\Omega+2)}\,u^2\right)\,
\bigl(1+o(1)\bigr),\qquad n\to\infty.
\label{eq:expansion-exp}
\end{equation}

\smallskip
\emph{(ii) Denominator.)} Set
\[
J(\lambda):=(1-\lambda)+\sqrt{\lambda/\Omega}\,R(z)
\]
where $R(z)$ is given by \eqref{eq:R(z)-def} with $z=2\sqrt{\lambda/\Omega}.$
At $\lambda=\lambda_0$, one checks that $J(\lambda_0)=0$. Differentiating (with
$s(\lambda):=\sqrt{1+4\lambda/\Omega}$),
\[
J'(\lambda)=-1+\frac{2/\Omega}{1+s}-\frac{(2\lambda/\Omega)\,\frac{2}{\Omega s}}{(1+s)^2},
\]
thereby
\begin{equation}
J'(\lambda_0)=-\,\frac{\Omega+1}{\Omega+2}.
\label{eq:Jprime}
\end{equation}
Therefore $J(\lambda)=J'(\lambda_0)(\lambda-\lambda_0)+o(\lambda-\lambda_0)$.

\medskip
Now, using
\[
(1+z_0^2)^{-1/4}=\left(1+\frac{4\lambda_0}{\Omega}\right)^{-1/4}
=\left(\frac{\Omega+2}{\Omega}\right)^{-1/2}
\]
and
\[
\frac{1+\sqrt{1+z_0^2}}{z_0}=\frac{\Omega+1}{\sqrt{\Omega+1}}=\sqrt{\Omega+1},
\]
together with $\sqrt{\lambda_0/\Omega}=\sqrt{(\Omega+1)}/\Omega$ and
$\nu=n/2$, we obtain from \eqref{eq:LHS-structure}
\begin{align}
\text{LHS}\sim&
-\,\frac{1}{\sqrt{\pi n}}\,
\frac{\sqrt{\lambda_0/\Omega}}{J(\lambda)}\,
\left(1+\frac{4\lambda_0}{\Omega}\right)^{-1/4}
\frac{1+\sqrt{1+z_0^2}}{z_0}\nonumber\\
&\times
\exp\!\left(-\,\frac{\Omega}{4(\Omega+2)}\,u^2\right),
\label{eq:LHS-near}
\end{align}
that is,
\begin{align}
\text{LHS}\sim &
-\,\frac{1}{\sqrt{\pi n}}\,
\frac{\displaystyle \frac{\sqrt{\Omega+1}}{\Omega}\cdot
\sqrt{\frac{\Omega}{\Omega+2}}\cdot \sqrt{\Omega+1}}{J(\lambda)}\nonumber\\
&\times \exp\!\left(-\,\frac{\Omega}{4(\Omega+2)}\,u^2\right)\nonumber\\
&
=
-\,\frac{1}{\sqrt{\pi n}}\,
\frac{\displaystyle \frac{\Omega+1}{\sqrt{\Omega(\Omega+2)}}}{J(\lambda)}\;
e^{-\,\frac{\Omega}{4(\Omega+2)}u^2}.
\label{eq:LHS-simplified}
\end{align}
Using $J(\lambda)\sim J'(\lambda_0)(\lambda-\lambda_0)=J'(\lambda_0)\,u/\sqrt{n}$ and
\eqref{eq:Jprime}, we finally get
\begin{equation}
\text{LHS}\sim
\frac{1}{\sqrt{\pi}}\sqrt{\frac{\Omega+2}{\Omega}}\;
\frac{e^{-\frac{\Omega}{4(\Omega+2)}u^2}}{u}\,.
\label{eq:LHS-final-shape}
\end{equation}

Now, equating \eqref{eq:LHS-final-shape} to $\epsilon$ yields
\begin{equation}
\epsilon\;\approx\;\frac{1}{\sqrt{\pi}}\sqrt{\frac{\Omega+2}{\Omega}}\;
\frac{\exp\!\left(-\dfrac{\Omega}{4(\Omega+2)}u^2\right)}{u}
\label{eq:eps-balance}
\end{equation}
with $u=\sqrt{n}(\lambda-\lambda_0)$.
Solving \eqref{eq:eps-balance} leads to
 the two-term large-$n$ approximation of $\lambda$ as
\begin{align}
\lambda \;&\approx\; 
1+\frac{1}{\Omega}\;+\;\frac{1}{\sqrt{n}}\,
\sqrt{\frac{2(\Omega+2)}{\Omega}\;
W\!\left(\frac{1}{2\pi\epsilon^2}\right)},
\quad (n\to\infty),
\label{eq:final}
\end{align}
where $W(\cdot)$ is the principal branch of the Lambert $W$-function (so that $\lambda>\lambda_0$ and the left-hand side of
\eqref{eq:main} is positive).
\medskip

\section{Simplification of The Converse Bound For Large Blocklength}\label{sec-app-R}

Consider
\begin{equation}\label{eq:def-En}
E_n\;=\;-\frac{1}{n}\log_2\!\Bigg(
P_{0,{\rm FA}}\left(\frac{n\lambda}{1+\Omega}\right) 
\Bigg),
\end{equation}
where $P_{0,{\rm FA}}\left(\frac{n\lambda}{1+\Omega}\right) $ is given by \eqref{eq:F_FA_bound1a}.
Using the uniform Debye (large $\nu$) approximation for fixed $z>0$ given by \eqref{eq:debye} and using \eqref{eq:ratio}, we can write the term in the logarithm as 
\begin{align}\label{eq:main-eta}
E_n
&= \frac{1}{\ln 2}\Bigg[
\frac{2+\lambda+1/\Omega+\Omega}{2(1+\Omega)}
-\frac{1}{4}\ln\!\Big(\frac{\lambda\Omega}{(1+\Omega)^2}\Big)\nonumber\\
&-\frac{1}{2}\,\eta\!\Big(2\sqrt{\tfrac{\lambda}{\Omega}}\Big)
\Bigg] + o(1),
\qquad n\to\infty,
\end{align}
with $\eta(z)$ given by \eqref{eq:debye1}. Inserting $\eta(z)$ from  \eqref{eq:debye1} for $z=2\sqrt{\lambda/\Omega}$ and simplifying,
we obtain the final form as
\begin{align}\label{eq:final-closed}
E_n
&= \frac{1}{2}\log_2\!\Big(1+\Omega \Big)+
\frac{1}{2\ln 2}\frac{2+\lambda+1/\Omega+\Omega}{1+\Omega}\nonumber\\
&
-\frac{1}{2\ln 2}\sqrt{1+\frac{4\lambda}{\Omega}}-\frac{1}{2}\log_2\!\left(\frac{2\lambda}{ 1+\sqrt{1+\frac{4\lambda}{\Omega}}}\right)
 + o(1).
\end{align}

\bibliographystyle{IEEEtran}
\bibliography{Citations}

\end{document}